\documentclass[iicol,pdflatex,sn-mathphys-ay]{sn-jnl}% Math and Physical Sciences Author Year Reference Style
\usepackage{graphicx}%
\usepackage{subcaption}
\usepackage{multirow}%
\usepackage{amsmath,amssymb,amsfonts}%
\usepackage{amsthm}%
\usepackage{mathrsfs}%
\usepackage[title]{appendix}%
\usepackage{xcolor}%
\usepackage{textcomp}%
\usepackage{manyfoot}%
\usepackage{booktabs}%
\usepackage{algorithm}%
\usepackage{algorithmicx}%
\usepackage{algpseudocode}%
\usepackage{listings}%
\usepackage{tikz}
\usepackage{pgfplots}
\pgfplotsset{compat=1.18}
\usepackage{stfloats}
\newcommand{\bs}{\boldsymbol}

\newcommand{\Ne}{{\text{N}_\text{e}}}
\newcommand{\Ns}{{\text{N}_\text{s}}}
\newcommand{\Nsf}{{\text{N}_\text{s}^f}}
\newcommand{\Nsg}{{\text{N}_\text{s}^g}}
\newcommand{\Nf}{{\text{N}_\text{f}}}
\newcommand{\Nu}{{\text{N}_\text{u}}}

\newcommand{\sens}[2]{\frac{\partial{#1}}{\partial{#2}}}
\newcommand{\scnd}[2]{\frac{\partial^2{#1}}{\partial{#2}^2}}
\newcommand{\scndcross}[3]{\frac{\partial^2{#1}}{\partial{#2}\,\partial{#3}}}
\newcommand{\drho}[1]{\frac{\partial{#1}}{\partial \rho_e}}
\newcommand{\total}[2]{\frac{\text{d}\,{#1}}{\text{d}\,{#2}}}
\newcommand{\dt}[1]{\frac{\partial{#1}}{\partial t}}
\newcommand{\dtt}[1]{\frac{\partial^2{#1}}{\partial t^2}}

\newcommand{\ba}{{\bs a}}
\newcommand{\bg}{{\bs g}}
\newcommand{\boldf}{{\bs f}}
\newcommand{\bolds}{{\bs s}}
\newcommand{\bn}{{\bs n}}
\newcommand{\bu}{{\bs u}}
\newcommand{\bx}{{\bs x}}
\newcommand{\bz}{{\bs z}}
\newcommand{\bA}{{\bs A}}
\newcommand{\bB}{{\bs B}}
\newcommand{\bG}{{\bs G}}
\newcommand{\bH}{{\bs H}}
\newcommand{\bK}{{\bs K}}
\newcommand{\bW}{{\bs W}}
\newcommand{\brho}{{\bs \rho}}
\newcommand{\blmbd}{{\bs \lambda}}
\newcommand{\rhomin}{\rho_{\min}}

\newcommand{\figref}[1]{Fig.~\ref{#1}}
\newcommand{\secref}[1]{Sec.~\ref{#1}}
\newcommand{\eqnref}[1]{(\ref{#1})}

\theoremstyle{thmstyleone}%
\theoremstyle{thmstyletwo}%

\theoremstyle{thmstylethree}%

\newcommand{\R}{\mathbb{R}}
\newcommand{\norm}[1]{\lVert#1\rVert}
\newcommand{\abs}[1]{\lvert#1\rvert}

\begin{document}

\title{Aligning feature-mapping designs to given density fields, first and second order}

\author*[1]{\fnm{Fabian} \sur{Wein}}\email{fabian.wein@fau.de}

\author{\fnm{Patrick} \sur{Jung}}\email{patrick.jung@fau.de}

\author[2]{\fnm{Arash} \sur{Moradian}}\email{arash.moradian@fau.de}

\author[1]{\fnm{Michael} \sur{Stingl}}\email{michael.stingl@fau.de}

\affil*[1]{\orgdiv{Department of Mathematics},
  \orgname{Friedrich-Alexander-Universit\"at Erlangen-N\"urnberg},
  \orgaddress{\city{Erlangen}, \country{Germany}}}

\affil*[2]{\orgdiv{Institute for Multiscale Simulation},
  \orgname{Friedrich-Alexander-Universit\"at Erlangen-N\"urnberg},
  \orgaddress{\city{Erlangen}, \country{Germany}}}

\abstract{%
Feature-mapping methods represent structural designs by explicit geometric
primitives on fixed analysis grids, an interpretable and parametric
alternative to the density fields of density-based topology optimization
(SIMP). In contrast to SIMP, however, feature-mapping optimization is sensitive to the initial design.

We present a gradient-based approach to align a feature-mapping configuration to a given (pseudo) density field, e.g., to initialize a subsequent feature-mapping optimization. The staged approach is based on a least-squares tracking formulation, preceded by a variant that only rewards alignment and does not penalize mismatch. Without an underlying finite element simulation, iterations are cheap, but isolated features receive little sensitivity information. As a remedy, we propose an asymmetric transition function based on automatically parametrized B\'ezier curves. We demonstrate the approach for a 2D cantilever and the five-bar design; for the latter, we show optional feature minimization based on the feature scaling variable of the geometry projection method.

The approach works well with first-order optimizers. The absence of a state problem and the small number of variables, however, also make a second-order formulation attractive. To the best of the authors' knowledge, we present the first Hessian formulation for feature mapping and report the behavior of first- and second-order optimizers on our benchmark problems. For completeness, we also give the exact Hessian of the state-based compliance, which requires one additional solution of the FEM system per feature variable.}


%%\pacs[MSC Classification]{49M05, 65K10, 74P05}

\maketitle

% ==========================================================================
\section{Introduction}\label{sec:introduction}
% ==========================================================================

Density-based topology optimization, and in particular the SIMP method,
is a very powerful and efficient approach to obtain highly effective structural layouts on fixed finite element grids. However, there are also use cases where the feature-mapping paradigm is preferred, when the design is composed of a small number of high-level geometric primitives like bars. Feature-mapping methods are based on the same pseudo-density $\brho \in \mathbb{R}^\Ne$ representation as SIMP, both for solving the state problem and for computing sensitivities of cost and constraint functions with respect to the $\Ne$ element-wise pseudo density variables. Feature mapping adds an additional layer, where the pseudo density field $\brho(\bolds)$ is constructed in a differentiable manner by the shape parameters $\bolds \in \mathbb{R}^\Ns$ for the features with $\Ns \ll \Ne$. See \cite{wein2020review} for an extensive analysis and review of feature-mapping methods and \secref{sec:features} for a brief introduction.

The efficiency of SIMP-based approaches stems from the gradient information -- cheap to compute, yet information-rich: with a value at every element, solid material and voids can emerge potentially anywhere at any iteration.
Feature mapping combines aspects from density-based optimization with shape optimization. Via chain rule the SIMP gradient field is restricted to the area around the feature boundary (\textit{transition zone}) and much of the SIMP gradient information is discarded. 
A transition zone that is too large blurs the design and may prevent a solid (fully dense) solution for thin features, while one that is too small reduces the information captured by the shape sensitivity of the feature parameters. This is one of the key drawbacks of feature-mapping methods, the often comparatively poor convergence and high sensitivity to the initial design for feature-mapping problems (number of features and their distribution in the design space).

Already in \citet{wein2020review} the tendency of feature mapping to become trapped in local minima is demonstrated. \citet{Zelickman2022} randomly generate 100 initial designs for the optimization of plate supports with a feature-mapping parametrization, of which only a quarter converge to the layout they identify as the global optimum; a combination of remedies restores reliability. Random-start studies with the tracking formulation are reported in \citep{Jung2026arXiv}.

The remedies proposed against this sensitivity fall into three groups. The first group leaves the initial design arbitrary and hardens the optimization instead, by continuation of the projection parameters and damping of the derivatives \citep{Zelickman2022}, by tunneling from one local minimum to a better one \citep{ZhangNorato2018}, or by statistical and surrogate-based optimization due to the small number of design variables (see \citep{wein2020review} for a survey). The second group makes the initial layout less relevant by changing the set of features during the optimization, either by inserting and deleting components \citep{Cui2022,Li2026CAID} or by constraining a smooth count of the effective features \citep{Zhang2017complexity}. The third group computes an initial layout from a density field: \citet{Weiss2020} threshold a SIMP result to a black-and-white design, skeletonize it to obtain the medial axis and connect the skeleton nodes by straight bars; related geometric pipelines, partly combined with sparse optimization, are used in \citep{Lian2020,Ling2024}, and learned variants recognize trusses or segment instances in the density image \citep{Gamache2018,Rochefort2023}.

Our work belongs to this third group. One motivation is the generation of (initial) feature-mapping designs from a given pseudo density field, obtained by a classical SIMP problem. In contrast to the approaches above (thresholding, skeletons, curve fitting, training), we replace the geometric pipeline by a sensitivity-based optimization problem in the feature parameters themselves. The formulation is furthermore not restricted to a particular feature type, whereas skeleton-based pipelines are tied to bar-like features.

Another motivation for the translation of a SIMP result to special high level features is to obtain a parametric CAD model for manual post-processing or possibly even the re-parametrization of feature-mapping result based on its pseudo density field,
similar to \cite{Shannon2023}, where post-processing pipelines convert feature-mapping results into parametric CAD models through feature removal, merging and snapping. One might also want to obtain an explicit and precise design representation of a SIMP result, e.g., if the problem is subject to boundary effects \citep{Wadbro:2026:BoundaryReview}. For an overview of the geometric post-processing of topology optimization results we refer to \citep{Subedi2020}, and to \citep{Koemma2026} for a fully automated interpretation of a topology optimization result into a parametric frame structure made of standard profiles. It is this second motivation that our numerical results address directly; a demonstration that a subsequent feature-mapping optimization started from a tracked design converges better than from a generic one is beyond the scope of this work.

Our basic idea is to match a given pseudo density field $\brho^*$ by the field $\brho(\bolds)$ generated by shape parameters $\bolds$. One can directly formulate this as a differentiable optimization problem
\begin{align}
\min_{\bolds} J_{\text{track}}(\bolds) = \norm{\brho(\bolds) - \brho^*}^2.
\label{eqn:tracking_rho}
\end{align}
This already works in principle, but as we show in this work, a staged approach with a modified objective function (\textit{reward function}) as a first step is beneficial. 

The mentioned localization of the shape sensitivity to the transition zone is even more severe in the present design tracking problem. With a physics based feature-mapping optimization problem, the density gradient $\drho{J}$ is a global field, i.e. load and support typically have a global influence and hence generally still impact within the transition zone. For plain design tracking, the shape gradient is perfectly blind outside the transition zone. The naive remedy, to enlarge the transition zone, fails for the standard symmetric transition function and small features. We propose an asymmetric transition extension that generates far-field sensitivity. 

Since the tracking objective depends on the feature parameters only through the explicit density map $\brho(\bolds)$, exact second-order sensitivities are accessible at low cost. For a physics-based feature-mapping problem, the Hessian would require the full $\frac{\partial^2 J}{\partial \rho_e \partial \rho_i} \in \mathbb{R}^{\Ne \times \Ne}$ operator, which is prohibitively expensive when obtained exactly. Together with the small number of shape parameters $\Ns \ll \Ne$, this makes Newton-type solvers with exact Hessians practical for design tracking -- provided that every link of the density map is twice differentiable, which holds for our proposed asymmetric transition function. For completeness we also give the formulation for compliance, where additional solutions of the FEM system are required.

This paper builds on the Master's thesis \cite{Jung2026arXiv}, which contains more detailed derivations of the sensitivity expressions, extended parameter studies and an alternative heuristic consolidation approach. However, the present work extends it significantly with respect to the asymmetric transition function and the Hessian formulations.
The main contributions of this work are
\begin{enumerate}
\item a least-squares tracking formulation that fits a feature-mapping design to
  a given pseudo density field, together with a reward variant which only
  rewards alignment and does not penalize mismatch,
\item a B\'ezier-curve-based asymmetric transition function with an automatic
  parametrization, which provides sensitivity information beyond the symmetric
  transition zone,
\item exact analytic Hessians for the state-less tracking problem, and the
  corresponding formulation for the state-based compliance function (to the
  best of the authors' knowledge the first second-order formulation of feature
  mapping),
\item an optional consolidation stage which removes unnecessary features by
  minimizing the sum of the fading variables under a bound on the tracking
  error.
\end{enumerate}

The remainder of the paper is organized as follows.
Section~\ref{sec:features} presents the feature model, the B\'ezier-based asymmetric transition function and the analytic
sensitivity chain.
Section~\ref{sec:opt} formulates the tracking and reward objectives.
Section~\ref{sec:second_order} derives the exact Hessian of the state-less problem
and its generalization to state-based functions.
Section~\ref{sec:geometry_variable} describes an optional feature consolidation approach.
Section~\ref{sec:stages} combines these ingredients into the staged strategy.
Section~\ref{sec:results} reports numerical results, and
Section~\ref{sec:conclusion} concludes.

% ==========================================================================
\section{Feature mapping}
\label{sec:features}
% ==========================================================================
In this section, we give a brief introduction to a standard 2D feature-mapping model based on capsule-shaped bar features. The 3D case is analogous without loss of generality. We refer to \cite{wein2020review} for a comprehensive introduction to the topic. The term \textit{feature mapping} subsumes a class of related methods, e.g., \textit{method of moving morphable components} (MMC) presented in \cite{Zhang:2016:MMC} (for variants see the MMC review paper \cite{Li:2024:MMCReview}), and \textit{geometry projection} (GP), \cite{Norato2015}, to name the most prominent ones. Two of the present authors used a spline-based shape description as a feature model \citep{WeinStingl2018}. Within the specific methods, terminology and realization details vary from the presentation here, but the principles are generally the same. We start with the standard symmetric polynomial transition function and later introduce our B\'ezier-curve-based asymmetric version. \figref{fig:sym_1d} summarizes the complete mapping in a 1D sketch: the signed distance field $d$ of the feature, the transition function $H$ smearing the boundary over a zone of half-width $a$, and the resulting pseudo-density $\rho$. The following subsections detail each of these steps.

\subsection{Exemplary single feature mapping}
\label{sec:geometry}
\begin{figure}
  \centering
  \includegraphics[width=0.39\linewidth]{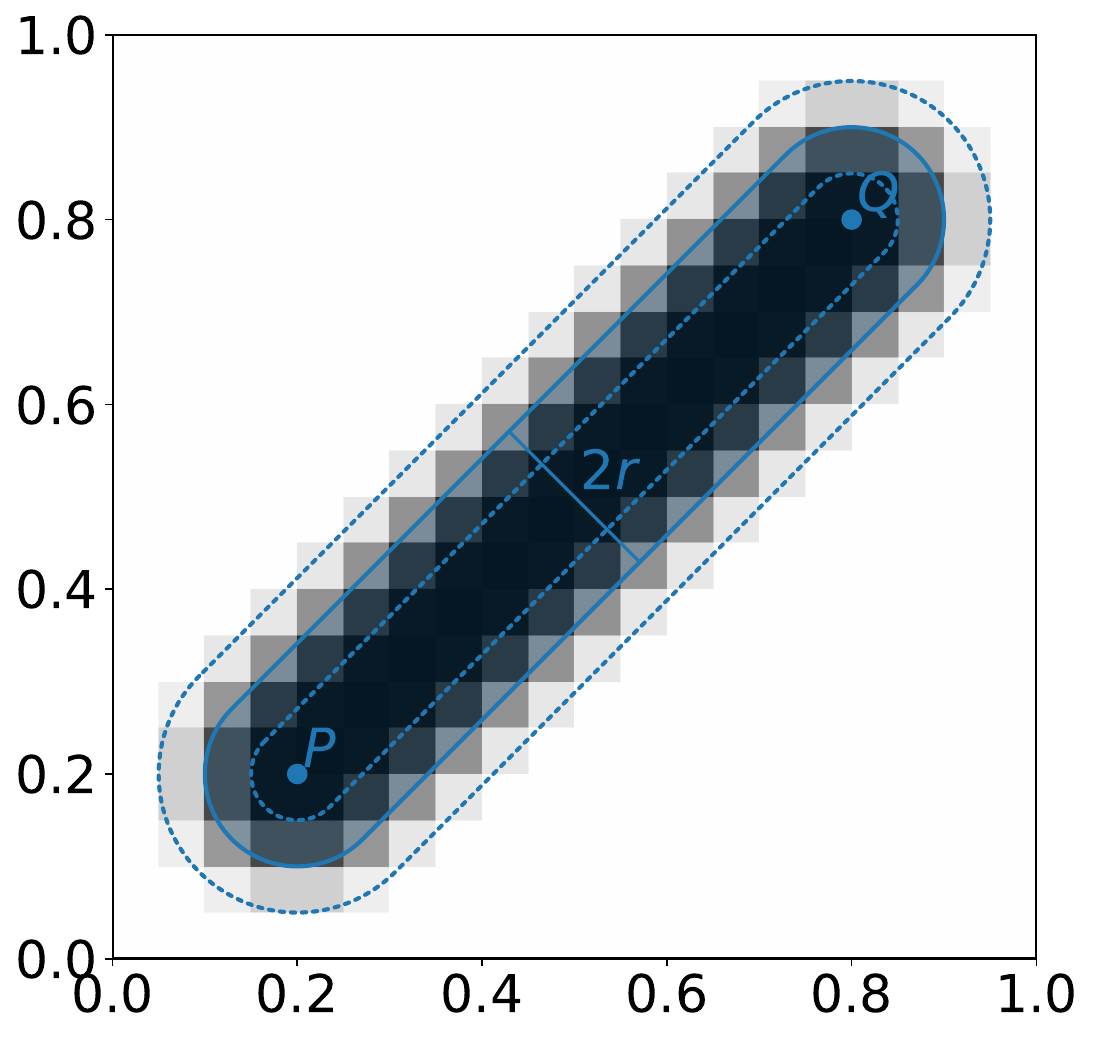}
  \includegraphics[width=0.59\linewidth]{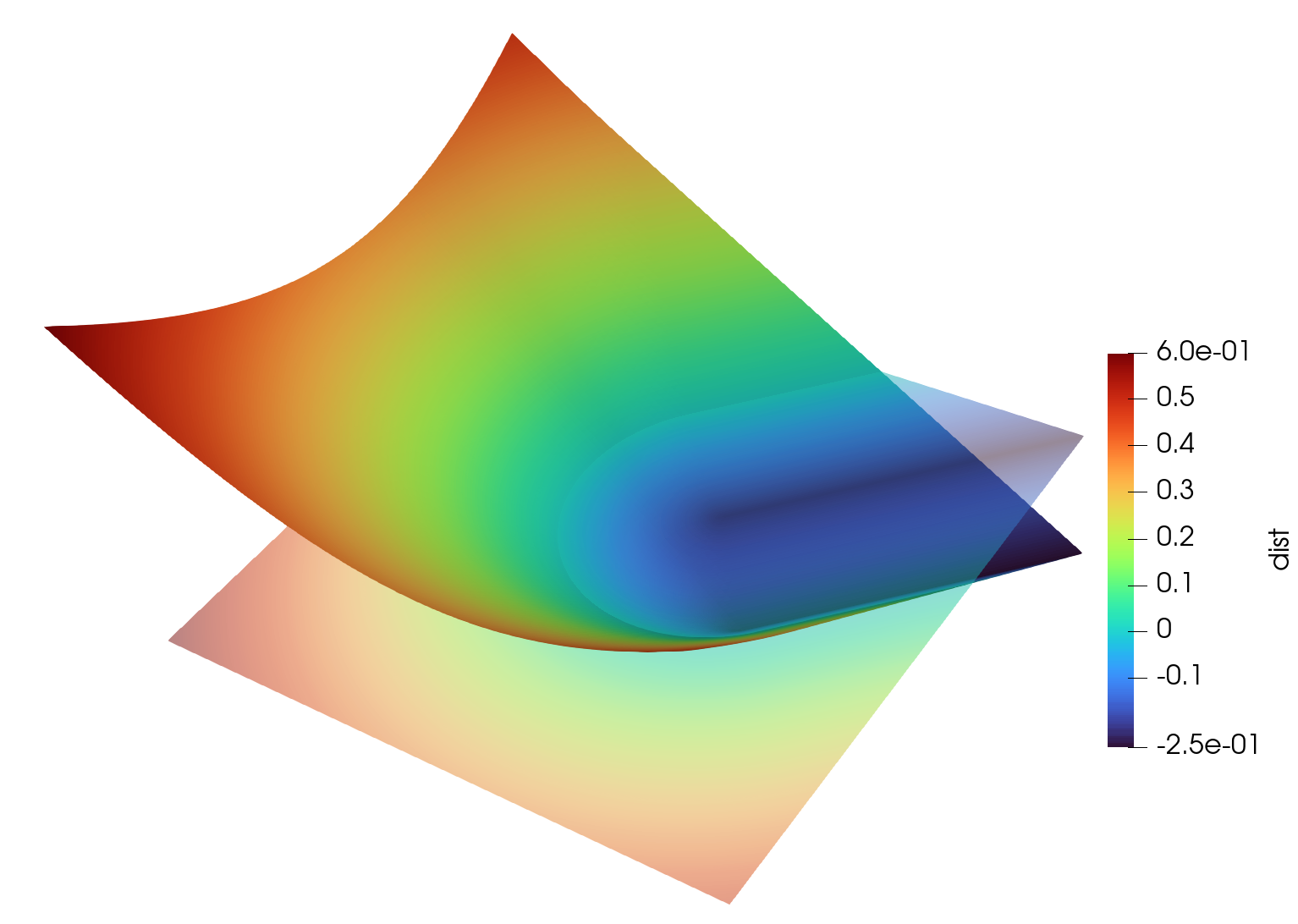}
  \caption{Left: a single bar feature and the resulting pseudo-density field. Right: signed distance field visualized as a height map.}
  \label{fig:single_density}
\end{figure}

The geometric primitive used throughout this work is a capsule-shaped bar with endpoints $P=(p_x,p_y)$, $Q=(q_x,q_y)$ and profile radius $r$, i.e.\ $r$ is the half-width of the capsule and $2r$ its full width. For the single feature case, this defines the shape variables
\begin{equation}
 \bolds = (p_x,p_y,q_x,q_y,r)^\top \in\R^{5},
 \label{eqn:shape_vars}
\end{equation}
see Fig.~\ref{fig:single_density}. In structural optimization this is an obvious model, e.g., employed in \cite{Norato2015}. The other common bar representation is by a rectangle or trapezoid; however, its convex corners come with minor differentiability issues and, unlike the $C^1$ capsule boundary, do not admit level sets parallel to the boundary everywhere.

\begin{figure*}[ht]
\centering
\begin{tikzpicture}[>=latex, font=\small]

% Layout: object center xc=7, radius r=3 -> boundary at x=4 and x=10
%         transition half-width a=1.5
%         d=-a (inner edge): x=5.5 and x=8.5
%         d=+a (outer edge): x=2.5 and x=11.5

% --- Axis ---
\draw[black, line width=0.6pt] (-0.3,0) -- (13.3,0);
\foreach \xi in {0,1,...,13}
  \draw[black] (\xi,0.07)--(\xi,-0.07);
\node[black, below=7pt] at (1,0) {1D element spacing};

% --- Signed distance (red), d(x) = |x-7| - 3, scaled by 0.30 ---
\draw[red!80!black, line width=1.5pt]
  (0,  {(7-0-3)*0.30}) --
  (7,  {-3*0.30})       --
  (13, {(13-7-3)*0.30});
\node[red!80!black, above] at (1.5,{(7-1.5-3)*0.30+0.05})
  {$d(\bx;\bolds)$};

% --- Transition H (green, cubic smoothstep) ---
% H = 0 far outside
\draw[green!55!black, line width=1.5pt] (-0.3,0) -- (2.5,0);
% Left transition: t=(5.5-x)/3, H=2t^3-3t^2+1
\draw[green!55!black, line width=1.5pt]
  plot[domain=2.5:5.5, samples=60, variable=\x]
  ({\x},{2*((5.5-\x)/3)^3 - 3*((5.5-\x)/3)^2 + 1});
% Flat H=1 inside
\draw[green!55!black, line width=1.5pt] (5.5,1) -- (8.5,1);
% Right transition: t=(x-8.5)/3, H=2t^3-3t^2+1
\draw[green!55!black, line width=1.5pt]
  plot[domain=8.5:11.5, samples=60, variable=\x]
  ({\x},{2*((\x-8.5)/3)^3 - 3*((\x-8.5)/3)^2 + 1});
% H = 0 far outside
\draw[green!55!black, line width=1.5pt] (11.5,0) -- (13.3,0);
\node[green!55!black, above] at (7,1.08) {$H(d(\bx;\bolds))$};

% --- Object (blue, d<=0, i.e. x in [4,10]) ---
\draw[blue!70!black, line width=2.5pt] (4,0) -- (10,0);
\node[blue!70!black, above=2pt] at (7,0) {feature($\bolds$)};

% --- 0.5 markers ---
\foreach \xb in {4, 10}{
  \fill[black] (\xb,0.5) circle (1.5pt);
  \draw[black, dashed, thin] (\xb,0)--(\xb,0.5);
}
\node[right, black] at (10.08,0.5) {$0.5$};

% --- Brackets (orange) on both sides ---
% Left boundary (x=4): outer bracket [2.5,4] and inner bracket [4,5.5]
\foreach \xv in {2.5,4,5.5}
  \draw[orange!80!black, thin] (\xv,0)--(\xv,-0.38);
\draw[orange!80!black, <->, line width=0.9pt] (2.5,-0.35)--(4,-0.35);
\node[orange!80!black, below] at (3.25,-0.38) {$a$};
\draw[orange!80!black, <->, line width=0.9pt] (4,-0.35)--(5.5,-0.35);
\node[orange!80!black, below] at (4.75,-0.38) {$a$};
% Right boundary (x=10): inner bracket [8.5,10] and outer bracket [10,11.5]
\foreach \xv in {8.5,10,11.5}
  \draw[orange!80!black, thin] (\xv,0)--(\xv,-0.38);
\draw[orange!80!black, <->, line width=0.9pt] (8.5,-0.35)--(10,-0.35);
\node[orange!80!black, below] at (9.25,-0.38) {$a$};
\draw[orange!80!black, <->, line width=0.9pt] (10,-0.35)--(11.5,-0.35);
\node[orange!80!black, below] at (10.75,-0.38) {$a$};

\end{tikzpicture}
\caption{The feature-mapping principle in 1D: the object (blue) is described by its
  signed distance $d(\bx;\bolds)$ (red), which the symmetric cubic smoothstep
  transition function $H$ \eqnref{eqn:cubic_poly} (green) maps to a smooth pseudo-density. Only within the
  transition zone is the density neither 0 nor 1 and sensitive to the shape variables $\bolds$.}
\label{fig:sym_1d}
\end{figure*}
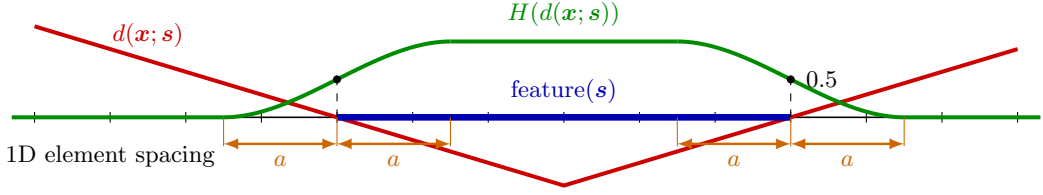

Crucial for feature mapping is the definition of a differentiable signed distance function $d(\bx; \bolds)$, which gives the shortest distance of a point $\bx$ to the feature boundary. For the present capsule feature, the distance is measured either to the nearer endpoint or, if the orthogonal projection of $\bx$ onto the axis falls within the segment, to the side. The relative position of this foot point,
\begin{equation}
  \beta(\bx;\bolds) = \frac{(\bx-P)\cdot(Q-P)}{\norm{Q-P}^2},
  \label{eqn:foot_point}
\end{equation}
selects the case,
\begin{equation}
  d(\bx;\bolds) =
  \begin{cases}
    \norm{\bx-P}-r & \text{if } \beta < 0 \\
    \norm{\bx-Q}-r & \text{if } \beta > 1 \\
    \abs{(\bx-Q)\cdot\bn}-r & \text{else,}
  \end{cases}
  \label{eqn:dist}
\end{equation}
where $\bn$ is perpendicular to the line segment from $P$ to $Q$ and normalized. As shown in Fig.~\ref{fig:single_density}, we obtain a continuous distance field, which is $C^1$ everywhere except for a kink at the line segment $\overline{PQ}$ (uncritical with $a<r$), and $C^2$ apart from the boundaries between the cap and side cases. There the gradient of $d$ stays continuous and only its second derivatives jump, as the level sets change from straight lines parallel to the axis within the side case to concentric circles within the caps. Outside the feature, the distance is positive, inside it is negative, and on the boundary it is zero, see Fig.~\ref{fig:sym_1d}. This is very similar to the level-set function $\phi(\bx)$ used in level-set methods (usually with opposite sign convention). The efficiency of the evaluation of the distance function is a key factor for the computational efficiency of feature-mapping methods. To apply different feature types, one just needs to provide an appropriate signed distance function.

The intention is to map the feature, given by its parameters $\bolds$, to a density field $\rho(\bx;\bolds) \in [0,1]$ with the idea of 1 (solid) inside the feature and 0 (void) outside. In structural optimization, a small positive lower bound $\rhomin$ is often used instead of zero. To obtain a differentiable mapping with $\sens{\rho}{s}$ available, we introduce a transition function $H:\R\to[0,1]$ that maps the signed distance to a (continuous) pseudo-density field. A common choice is a transition function symmetric around the feature boundary with a width of $2\,a$, see Fig.~\ref{fig:sym_1d}. In Sec.~\ref{sec:asym_transition} this will be extended to the asymmetric case. The requirements are $H(-a) = 1,\, H(0) = 0.5, \,H(a) = 0$ and $H'(-a)=H'(a)=0$. This is fulfilled by a cubic polynomial 
\begin{equation}
\label{eqn:cubic_poly}
H(d(\bx); a) :=
\begin{cases}
1 & \text{if } d(\bx) < -a \\
\frac{3}{4}\left( \frac{d(\bx)^3}{3 \, a^3} - \frac{d(\bx)}{a} \right) + \frac{1}{2} & \text{if }  |d(\bx)| \leq a \\
0 & \text{if } d(\bx) > a.
\end{cases}
\end{equation}
The final step for the mapping is to obtain piecewise constant pseudo density values $\rho_e \in [0,1]$ by integrating $H$,
% \begin{align*}
% \rho_e = \frac{1}{|\Omega_e|}\int_{\Omega_e} H(d(\bx;\bolds)) \, \intd x,
% \end{align*}
which is practically realized by numerical integration
\begin{equation}
\label{eqn:rho_e_num_int}
\rho_e = \frac{1}{N_\text{ip}} \sum_{i=1}^{N_\text{ip}} H(d(\bx_i;\bolds)).
\end{equation}
Many authors apply midpoint evaluation, i.e., $N_\text{ip}=1$, but as pointed out in \cite{wein2020review}, a sufficiently accurate numerical integration (e.g., $5 \times 5$ in 2D) indeed matters.

\subsection{First-order sensitivity analysis}
\label{sec:fm_sens}
The first-order analysis is standard, but we give it in some detail, as the intermediate derivatives are reused for the second-order analysis in \secref{sec:second-order_analysis}. This is the single-feature formulation, multi-feature aggregation is added in \secref{sec:aggregation}.

In full formulation, we have $\rho(\bx) = H(d(\bx;\bolds); a)$ and obtain the sensitivity via chain rule using \eqnref{eqn:rho_e_num_int}, \eqnref{eqn:cubic_poly} and \eqnref{eqn:dist}:
\begin{equation}
  \sens{\rho_e}{s_j} = \frac{1}{N_\text{ip}} \sum_{i=1}^{N_\text{ip}} \sens{H(d(\bx_i;\bolds); a)}{d} \sens{d(\bx_i;\bolds)}{s_j},
  \label{eqn:rho_sens}
\end{equation}
with 
\begin{equation}
  \sens{H(d(\bx); a)}{d} :=
  \begin{cases}
  0 & \text{if } d(\bx) < -a \\
  \frac{3}{4}\left( \frac{d(\bx)^2}{a^3} - \frac{1}{a} \right) & \text{if }  |d(\bx)| \leq a \\
  0 & \text{if } d(\bx) > a
  \end{cases}
  \label{eqn:H_sens}
\end{equation}
and $s_j$ being one of the shape parameters in \eqnref{eqn:shape_vars}. We label the three cases of \eqnref{eqn:dist} as
\begin{equation}
  d(\bx;\bolds) =
  \left\{
  \begin{array}{@{}l @{\;:\;} l @{}}
    \norm{\bx-P}-r & \text{cap}_P \\
    \norm{\bx-Q}-r & \text{cap}_Q \\
    \abs{(\bx-Q)\cdot\bn}-r & \text{side}.
  \end{array}
  \right.
  \label{eqn:dist_cases}
\end{equation}
We use the notation with $\text{cap}_P$, $\text{cap}_Q$ and $\text{side}$ as a marker for the derivatives of the corresponding part of the feature with the closest distance to $\bx$. The $\text{side}$ case expands to
\begin{equation}
\begin{split}
  (\bx-Q)\cdot\bn & = \frac{N(\bx; \bolds)}{D(\bolds)} = \frac{N}{\norm{P-Q}} \\
  & = \frac{(x-q_x)(p_y-q_y) + (y-q_y)(q_x-p_x)}{\sqrt{(p_x-q_x)^2+(p_y-q_y)^2}}.
\end{split}
  \label{eqn:dist_case_side}
\end{equation}
We need sensitivity expressions for each variable $s_j$ and each case in the distance function, e.g., for $s_0 = p_x$, we have, with the case selection of \eqnref{eqn:dist_cases},
\begin{equation}
  \sens{d(\bx)}{p_x} =
  \left\{
  \begin{array}{@{}l @{\;:\;} l @{}}
  \frac{p_x-x}{\norm{\bx-P}} & \text{cap}_P \\[1ex]
  0 & \text{cap}_Q \\[1ex]
  \frac{\operatorname{sgn}(N)\,\bigl(-(y-q_y)D^2-(p_x-q_x)N\bigr)}{D^3}& \text{side}.
  \end{array}
  \right.
  \label{eqn:sens_d_p_x}
\end{equation}
$\operatorname{sgn}(N)$ is non-differentiable at $N=0$. As the $\text{side}$ case applies only for $0 \leq \beta \leq 1$, this happens exclusively on the axis within the feature, where $a<r$ places the point outside the transition zone with $\sens{H}{d}=0$, hence it does not cause any issues. Furthermore, we need $D= \|P - Q\| > 0$. 
It is useful to store the case when evaluating the distance function. For comprehensive coverage of all cases, we refer to \cite{Jung2026arXiv}. This feature model is just an example of our tracking approach, and users can directly apply it to their favorite feature model by providing the corresponding distance function and its sensitivities.

\begin{figure*}[t]
  \centering
  \includegraphics[width=0.19\linewidth]{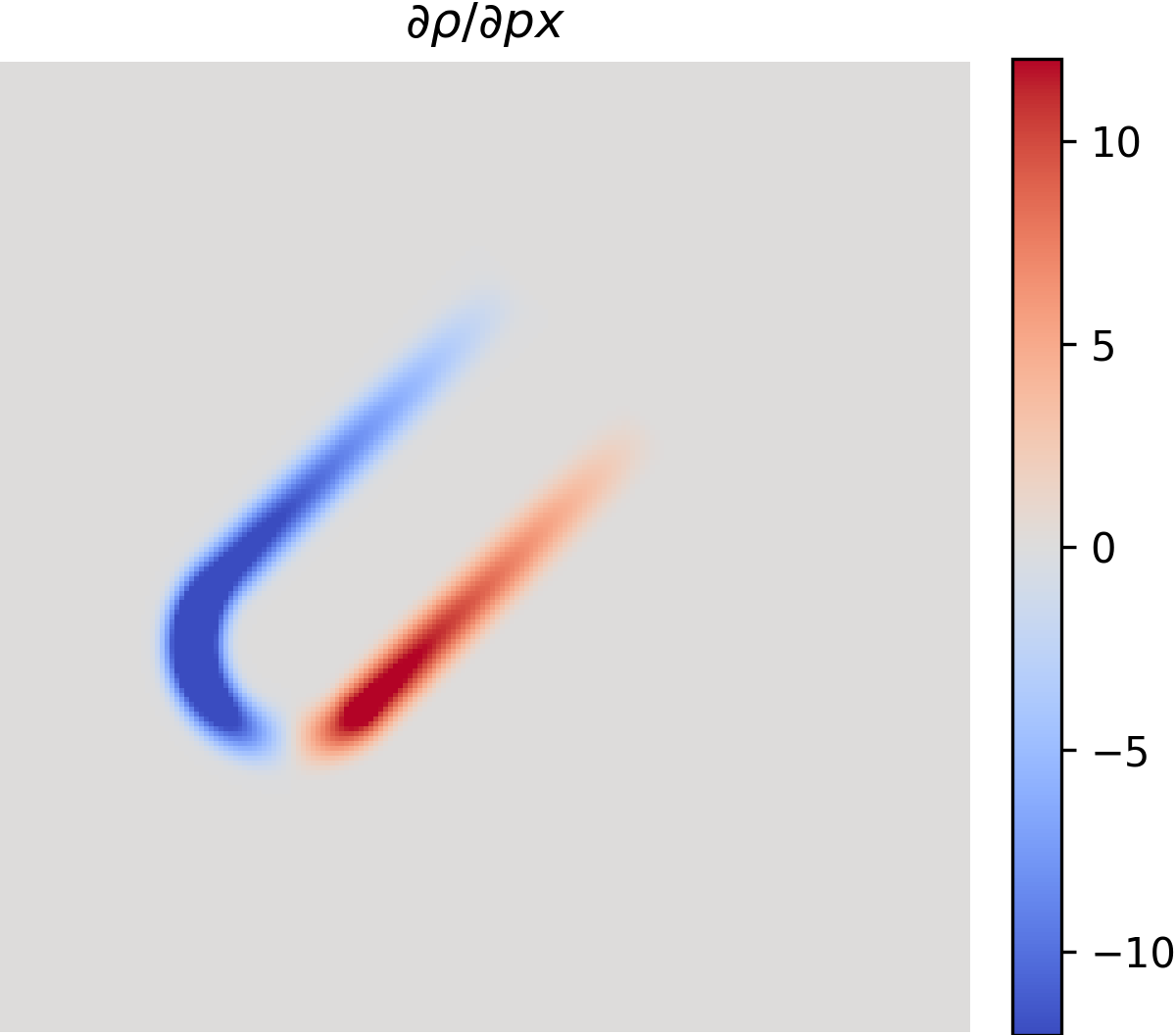}\hfill
  \includegraphics[width=0.19\linewidth]{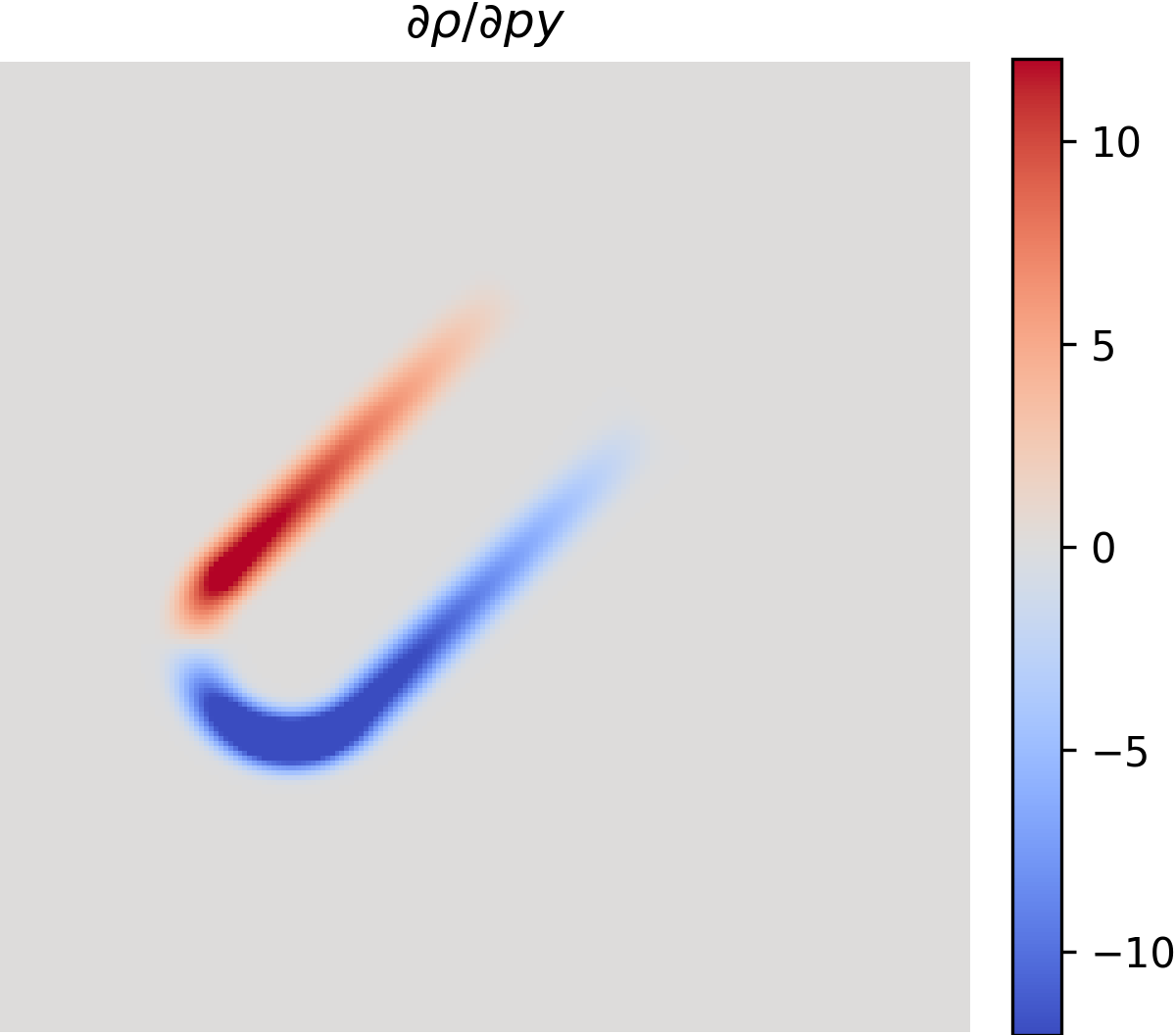}\hfill
  \includegraphics[width=0.19\linewidth]{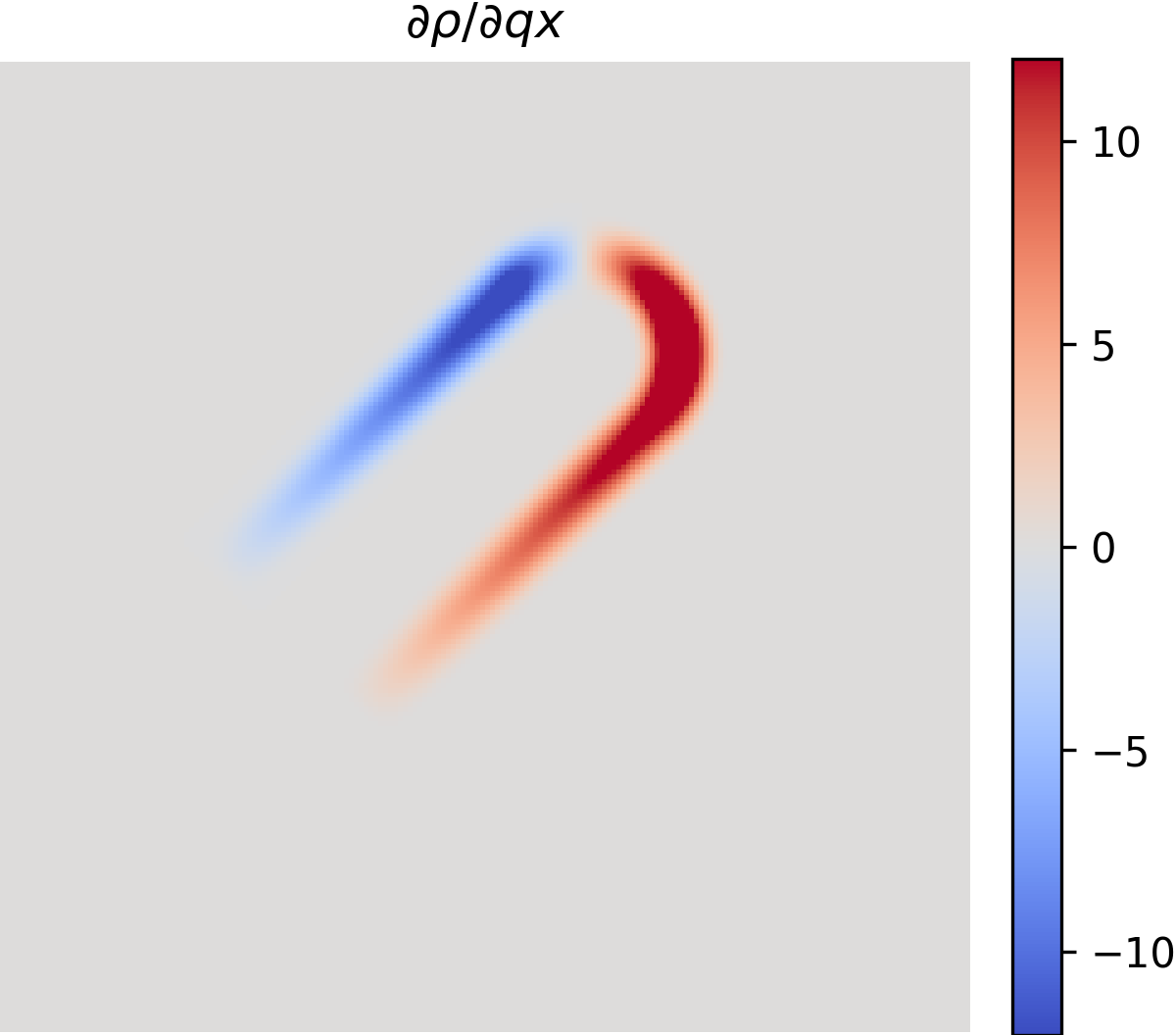}\hfill
  \includegraphics[width=0.19\linewidth]{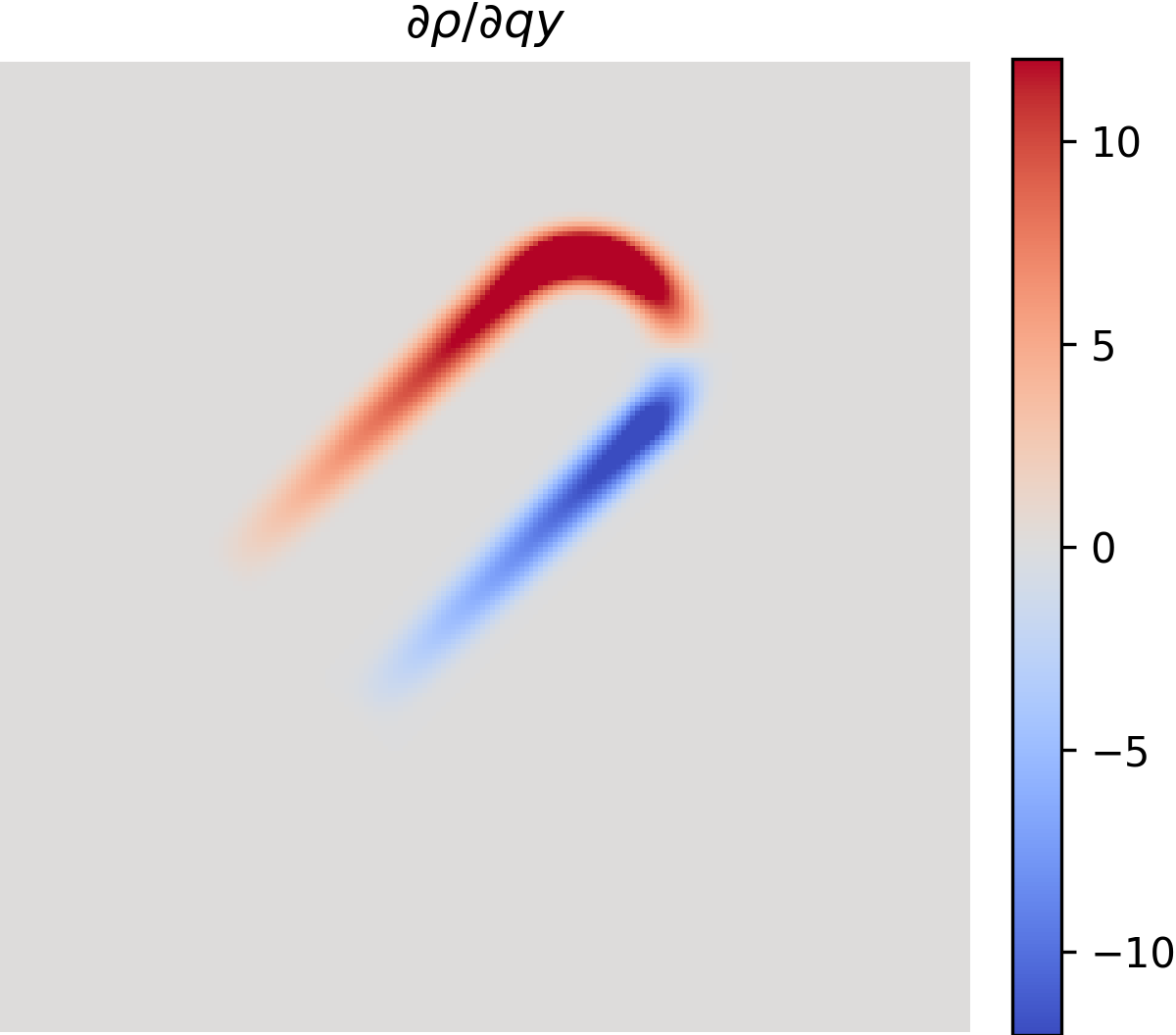}\hfill
  \includegraphics[width=0.19\linewidth]{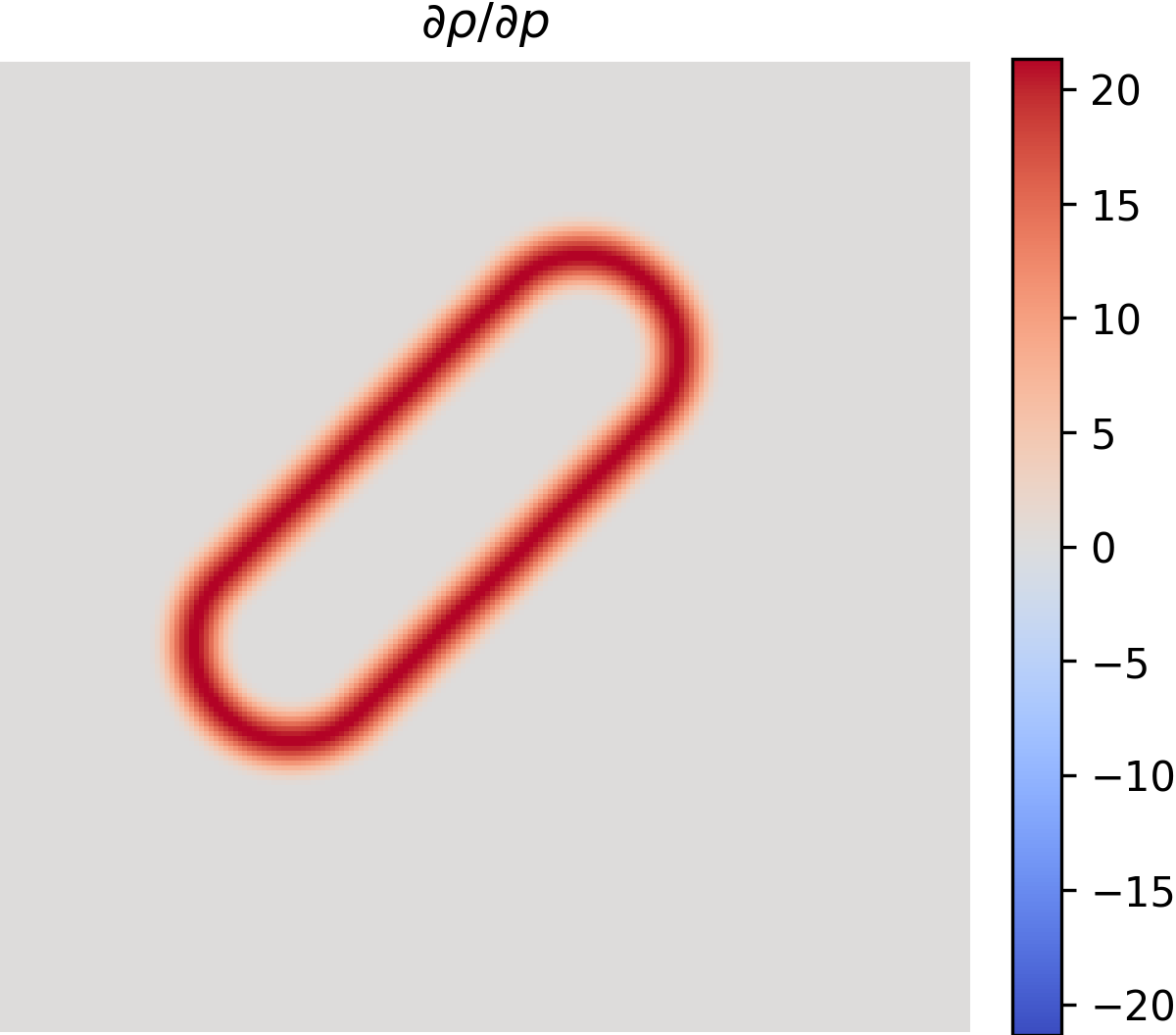}
  \caption{First-order sensitivities of a single feature, visualized for a high resolution of elements. The images show the fields $\sens{H}{s_j} = \sens{H}{d} \, \sens{d}{s_j}$, see \eqnref{eqn:H_sens} and, for $p_x$, \eqnref{eqn:sens_d_p_x}; the element sensitivity \eqnref{eqn:rho_sens} is their average over the integration points. From left to right: $s_j = p_x,\, p_y,\, q_x,\, q_y,\, r$.}
  \label{fig:sensitivities}
\end{figure*}

From \eqnref{eqn:H_sens} we see that the sensitivity of $\rho$ \eqnref{eqn:rho_sens} is zero outside the transition zone. A specific pseudo density $\rho_e$ therefore reacts to a slight variation of $\bolds$ only if at least one of its integration points falls into the transition zone around the feature boundary. All other elements are blind to the shape variables. See Fig.~\ref{fig:sensitivities} for a visualization of $\sens{H}{d} \sens{d}{s}$.

\subsection{Asymmetric transition function}
\label{sec:asym_transition}
With a symmetric transition zone \eqnref{eqn:cubic_poly}, we cannot simply make it arbitrarily large since we want and need solid inside, hence we need $a < r$. We therefore introduce an asymmetric transition zone with $a$ inside and $b$ outside. We present a $C^1$ version for the standard first-order optimization in topology optimization but also a $C^2$ version, to be motivated in Sec.~\ref{sec:second_order}. Our requirements are monotonicity, smoothness,
\begin{equation}
\begin{gathered}
  H(-a)=1,\, H(0)= 0.5,\, H(b) = 0, \\
  H'(-a)=H'(b)=0
\end{gathered}
  \label{eqn:asym_c1_req}
\end{equation}
and for $C^2$ also
\begin{equation}
  H''(-a)=H''(b)=0.
  \label{eqn:asym_c2_req}
\end{equation}
In the symmetric case $a=b$, polynomials suffice: the cubic \eqnref{eqn:cubic_poly} for $C^1$ and a quintic one for $C^2$. For $a \neq b$, it is much more difficult to fulfill all requirements, especially monotonicity and $C^1$ or $C^2$ smoothness everywhere. We propose a B\'ezier-curve-based approach, which can be parametrized to fulfill all requirements under certain conditions on the parametrization. A general formulation for a B\'ezier curve of degree $n$ using Bernstein polynomials with
\begin{equation*}
  B_i^n(t) := \binom{n}{i} (1-t)^{n-i} \, t^i, \quad i=0,\ldots,n,
\end{equation*}  
is given by 
\begin{equation}
  \bB(t) = \begin{pmatrix} b_x(t) \\ b_y(t) \end{pmatrix} = \sum_{i=0}^n B_i^n(t) \,\bW_i, \quad t \in [0,1],
\end{equation}
with the curve parameter $t$ and control points $\bW_i$. Interpreted as a transition function, we have $b_x(t) = d$ and $b_y(t) = H(d)$, hence $H(d) = b_y(b_x^{-1}(d))$ for $d \in [-a,b]$, under the condition that $b_x$ is strictly increasing on $t \in [0,1]$. In practice, one does not analytically invert $t(d)=b_x^{-1}(d)$ but searches for it. For fast numerics, we sample $H(d)$ and its derivatives and interpolate using look-up tables. For proper convergence, high accuracy is required.

\begin{figure}[htp]
  \centering
  \includegraphics[width=\linewidth]{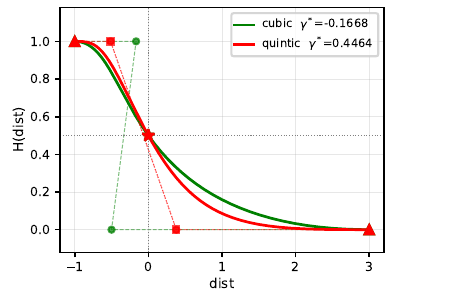} \\[1ex]
  \includegraphics[width=\linewidth]{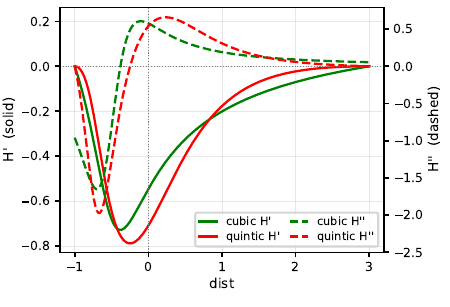}
  \caption{Top: asymmetric transition function $H(d)$ based on cubic and quintic B\'ezier curves with $a=1$ and $b=3$ and $\gamma^*$ chosen to minimize $\max_{d}\left|H''(d)\right|$. Bottom: the corresponding first and second derivatives with respect to the distance $d$.}
  \label{fig:asym}  
\end{figure}

The derivative of a B\'ezier curve of degree $n$ is again a B\'ezier curve of degree $n-1$,
\begin{equation*}
  \dt{\bB(t)} = n \sum_{i=0}^{n-1} B_i^{n-1}(t) (\bW_{i+1} - \bW_i).
\end{equation*}
We are actually interested in
\begin{equation}
  H'(d) = \sens{H(d)}{d} = \dt{b_y(t)} \cdot \sens{t}{d} = \frac{\dt{b_y(t)}}{\dt{b_x(t)}}
  \label{eqn:H_grad}
\end{equation} and
\begin{equation}
  H''(d) = \scnd{H(d)}{d} = 
  \frac{\dtt{b_y}\,\dt{b_x} - \dt{b_y}\,\dtt{b_x}}{\left(\dt{b_x}\right)^{3}}.
  \label{eqn:H_scnd}
\end{equation}
Note the third power. To obtain $H'(-a)=H'(b)=0$ for a cubic B\'ezier curve, it is sufficient to set the $y$-components of the control points $\bW_0$ to $\bW_3$ to 1, 1, 0, 0, and for a quintic one to 1, 1, 1, 0, 0, 0 to also obtain $H''(-a)=H''(b)=0$. We need $H(0)=0.5$ and ask this for $t=0.5$, hence $b_y(t=0.5)=0.5$ and $b_x(t=0.5)=0$. The first is given already by the chosen $y$-components of the control points. The second is fulfilled by setting the $x$-components of the cubic control points to $-a, c-\gamma, c+\gamma, b$ with
\begin{equation*}
  c = \frac{a-b}{\sum_{i=1}^{n-1} \binom{n}{i}}  
\end{equation*}
which is $c_\text{cubic} = \frac{a-b}{6}$ for the cubic case and we are left with a single parameter $\gamma$ to control the shape of the transition function and make sure that $b_x$ is 
invertible. For the quintic case we have two more control points $\bW_4$ and $\bW_5$ and could control their $x$-components by an additional parameter, but we choose to double the inner control points to have a single parameter $\gamma$ with $c_\text{quintic} = \frac{a-b}{30}$ and the $x$-components of the control points $-a, c-\gamma, c-\gamma, c+\gamma, c+\gamma, b$. See Fig.~\ref{fig:asym} for example curves.

\begin{figure}
  \centering
  \includegraphics[width=\linewidth]{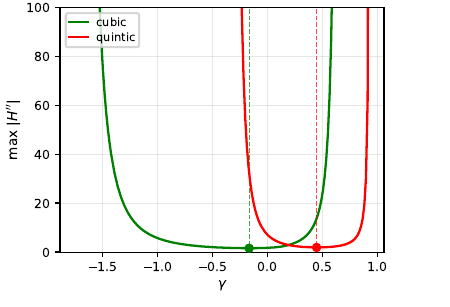} \\[1ex]
  \includegraphics[width=\linewidth]{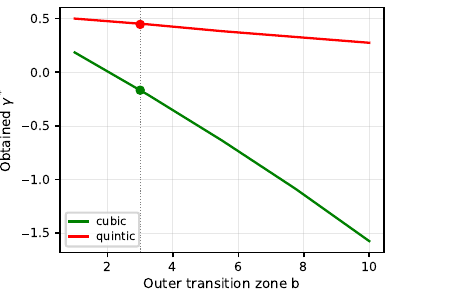}
  \caption{Obtaining the optimal parameter $\gamma^*$ for $a=1$ and $b=3$ as found by \eqnref{eqn:gamma_star} and a plot of $\gamma^*$ for a range of $b$.}
  \label{fig:asym_gamma}
\end{figure}

The parameters $a$ and $b$ are user choices to define the transition zone. A heuristic to define $\gamma$ for a given pair $a,b$ is to make the largest second derivative in the transition zone as small as possible. We apply it separately to the cubic and the quintic case,
\begin{equation}
  \gamma^* = \operatorname*{arg\,min}_{\gamma}
           \max_{d \in [-a,\, b]} \left|H''(d;\,\gamma)\right|.
  \label{eqn:gamma_star}         
\end{equation}
This leaves only the question of the range to search for $\gamma^*$. The $x$-coordinates of the inner control points need to lie inside the outer control points, hence we have $c-\gamma > -a$ and $c+\gamma < b$. Within this range we check numerically for $\dt{b_x(t; a, b, \gamma)} > 0 \; \forall \, t \in [0,1]$, see Fig.~\ref{fig:asym_gamma}.
If one instead minimized $\max_{d\in[-a,b]}|H'(d)|$, the solution for $a=b$ 
would approach a straight line, with the inner control points collapsing onto 
the outer ones. In this degenerate case, we have a kink at the endpoints, 
causing a singularity in the formula for $H''(d)$, making the B\'ezier 
parametrization unsuitable for computing second-order derivatives.

The transition function becomes, as an alternative to the symmetric polynomial-based form \eqnref{eqn:cubic_poly}, the following, which is to be efficiently evaluated based on precalculated data
\begin{equation}
\label{eqn:H_bezier}
\begin{aligned}
& H(d(\bx); a, b, \gamma^*) := \\
& \qquad
\begin{cases}
1 & \text{if } d(\bx) < -a \\
b_y(b_x^{-1}(d(\bx))) & \text{if } -a \leq d(\bx) \leq b \\
0 & \text{if } d(\bx) > b.
\end{cases}
\end{aligned}
\end{equation}

Both $a$ and $b$ are naturally measured in units of the element size $h$, because
the mapping \eqnref{eqn:rho_e_num_int} resolves the transition zone by sampling the
elements it crosses. This bounds $a$ from below: a transition zone much narrower than
the spacing of the integration points is sampled by too few of them, and the mapped
density $\rho_e$ then changes abruptly as a feature moves. Together with $a<r$ from
above, we obtain the scale chain
\begin{equation*}
  h \;\lesssim\; a \;<\; r .
\end{equation*}
The outer zone $b$ carries no such restriction and is a pure modeling choice: it is
the distance over which a feature senses a target, hence the parameter to enlarge
when features are initialized far away from the target regions. All numerical
examples in this work use $h=1/60$ with $a=3h$, so that the full symmetric
transition zone $2\,a$ spans six elements, an asymmetric extension of $12h$
giving $b=15h$, and profile bounds $r_\mathrm{min}=3.6h$ and $r_\mathrm{max}=30h$;
see \secref{sec:results}.

We provide an interactive JavaScript tool to explore the parametrization of the B\'ezier-based transition function and its derivatives at \url{https://am-ko.mi.uni-erlangen.de/bezier.html}.

\subsection{Aggregation of multiple features}
\label{sec:aggregation}
In realistic designs, many geometric features coexist and must be represented jointly. There are two principal approaches: map-then-combine and combine-then-map, see \cite{wein2020review} for more details. We follow here the more common and easier-to-implement map-then-combine approach. We assume $\Nf$ features, each consisting of $\Nsf$ shape parameters, and for simplicity we assume here all features to be capsules with $\Nsf=5$, see \eqnref{eqn:shape_vars}. The per-feature vectors $\bolds^f \in \R^{\Nsf}$ are collected in the design vector $\bolds \in \R^{\Ns}$ of dimension $\Ns = \sum_{f=1}^{\Nf} \Nsf = 5\,\Nf$. We use the common \emph{$p$-norm} aggregation
\begin{equation}
  \widetilde{\rho}_e = \left(\sum_{f=1}^{\Nf} \left(\rho_e^f\right)^p\right)^{1/p},
  \label{eqn:pnorm}
\end{equation}
where $\widetilde{\rho}_e$ is the aggregated density for element $e$ and $\rho_e^f$ are the integrated feature densities for all features $f=1,\ldots,\Nf$, evaluated by \eqnref{eqn:rho_e_num_int}.

Overlapping of $N$ features results in maximal values
\begin{align}
  \max \widetilde{\rho}_e = N^{1/p}.
  \label{eqn:p_norm_max}
\end{align}
For $p=1$ we have the simple sum. E.g., 4 features and $p=4$ give a maximum value of $\approx 1.414$ and for $p=8$ still $\approx 1.189$. For common topology optimization, e.g., compliance minimization, overlapping of features can therefore give a slight benefit. In our later tracking functional, we can make use of overshooting to gain a repulsive effect. The gradient is straightforward
\begin{equation}
  \sens{\widetilde{\rho}_e}{\rho_e^f} = \left(\sum_{g=1}^{\Nf} (\rho_e^g)^p\right)^{\frac{1}{p}-1} (\rho_e^f)^{p-1} = \left( \frac{\rho_e^f}{\widetilde{\rho}_e} \right)^{p-1}.
  \label{eqn:grad_pnorm}
\end{equation}
When overshooting, or the implicit effect of it, becomes an issue for standard problems, other aggregation functions are used, some are listed in \cite{wein2020review}. Another approach is to fade features via additional variables $\alpha_f \in [0,1]$, see \cite{Norato2015}. We take up this variable in \secref{sec:geometry_variable} as the basis of our optional consolidation stage.
\section{Optimization problem}
\label{sec:opt}

\subsection{Objectives}\label{sec:objectives}
The primary objective is a least-squares functional penalizing deviations between a given target density $\brho^*$ and the aggregated element density $\widetilde{\brho}(\bolds)$
\begin{equation}
  J_\mathrm{track}(\bolds)
  = \sum_{e=1}^\Ne \left(\rho^*_e - \widetilde{\rho}_e(\bolds)\right)^2 .
  \label{eqn:tracking_obj}
\end{equation}
This functional penalizes both uncovered target regions and features outside the target density. A local minimum for the latter, when the target density cannot be reached, is to shrink features or move them (partially) outside the domain. 
We therefore use a reward-only functional for early exploration
\begin{equation}
  J_\mathrm{reward}(\bolds)
  = -\sum_{e=1}^\Ne \rho^*_e \widetilde{\rho}_e(\bolds),
  \label{eqn:reward_obj}
\end{equation}
which maximizes the overlap with the target region and sees nothing else; the minus sign performs this maximization in a minimization framework. Note that expanding the features across the whole domain is a trivial optimum. We therefore keep the profile $r$ fixed while the reward functional is active, see \secref{sec:stages}.

The objective derivatives follow directly,
\begin{equation}
  \sens{J_\mathrm{track}}{\widetilde{\rho}_e}
  = -2\left( \rho^*_e - \widetilde{\rho}_e(\bolds)\right) \quad\text{and}\quad
  \sens{J_\mathrm{reward}}{\widetilde{\rho}_e} = -\rho^*_e.
  \label{eqn:obj_grad}
\end{equation}
%and enter the gradient with respect to the design variables through the chain rule,
% \begin{equation}
%   \total{J}{s_j^f} = \sum_{e=1}^{\Ne}
%   \sens{J}{\widetilde{\rho}_e}\,\sens{\widetilde{\rho}_e(\bolds)}{s_j^f} .
%   \label{eqn:obj_grad_chain}
% \end{equation}
% This contraction of the SIMP gradient field with the density Jacobian is the essence of feature mapping; \eqnref{eqn:gradient_chain} expands the Jacobian into its layers.
Note that $\sens{J_\mathrm{reward}}{\widetilde{\rho}_e}$ vanishes wherever $\rho^*_e = 0$, and it is the asymmetric transition zone \eqnref{eqn:H_bezier} which lets a distant feature reach a target element at all.

\subsection{Bounds and constraints}\label{sec:bounds}
Standard in feature mapping are constraints on the geometry of the features. Restricting the length of a feature $f$ by
\begin{equation}
  l^f(\bolds) = \norm{Q^f-P^f} = \sqrt{(q_x^f-p_x^f)^2 + (q_y^f-p_y^f)^2}, 
  \label{eqn:length}
\end{equation}
prevents the degenerate case $P^f=Q^f$. A not too small bound also prevents singularities, as this length appears as $D$ in the denominator of the shape derivatives \eqnref{eqn:sens_d_p_x}.
In the standard case we have a rectangular design domain
\begin{equation}
  \Omega = [x_{\min},\, x_{\max}] \times [y_{\min},\, y_{\max}] \subset \mathbb{R}^2
\end{equation}
and when we restrict all $P^f, Q^f \in \Omega$, we make sure the features cannot fully leave the domain, but parts of a feature can reach beyond.

\subsection{Full formulation}
\label{sec:full_formulation}
A full problem formulation for either objective reads
\begin{equation}
\begin{aligned}
  \min_{\bolds\,\in\,\R^{\Ns}}
     & J(\bolds), && J \in \{J_\mathrm{reward},\, J_\mathrm{track}\}\\
  \text{s.t.}\quad
      l^f(\bolds) &= \norm{Q^f-P^f} \ge \ell_\mathrm{min} && \forall\, 1\leq f \leq \Nf,\\
      P^f,\, Q^f &\in \Omega && \forall\, 1\leq f \leq \Nf,\\
     r^f &\in [\,r_\mathrm{min},\, r_\mathrm{max}\,] && \forall\, 1\leq f \leq \Nf.
\end{aligned}
\label{eqn:full_formulation}
\end{equation}
We require $r_\mathrm{min} > a$, such that each feature always retains a small solid core and the mapping remains differentiable. The full mapping chains the distance \eqnref{eqn:dist}, the transition function \eqnref{eqn:H_bezier}, the integration \eqnref{eqn:rho_e_num_int} and the aggregation \eqnref{eqn:pnorm} into one of the objectives \eqnref{eqn:reward_obj}, \eqnref{eqn:tracking_obj}:
\begin{equation}
\begin{aligned}
  \bolds^f
  &\xrightarrow{\ \text{distance}\ } d(\bx;\bolds^f)
   \xrightarrow{\ H\ } H(d;a,b,\gamma^*)\\
  &\xrightarrow{\ \text{integrate}\ }
     \rho_e^f = \frac{1}{N_\mathrm{ip}}\sum_{i=1}^{N_\mathrm{ip}} H\!\left(d(\bx_i;\bolds^f)\right)\\
  &\xrightarrow{\ \text{aggregate}\ }
     \widetilde{\rho}_e = \Big(\sum_{g=1}^{\Nf}(\rho_e^g)^p\Big)^{1/p}\\
  &\xrightarrow{\ \text{function}\ }
     J = \sum_{e=1}^{\Ne}\big(\rho^*_e-\widetilde{\rho}_e\big)^2 .
\end{aligned}
\label{eqn:forward_map}
\end{equation}
\noindent The full gradient follows directly by the chain rule
\begin{equation}
\begin{split}
  \total{J}{s_j^f} = \sum_{e=1}^{\Ne}
  &\underbrace{\sens{J}{\widetilde{\rho}_e}}_{\text{function}}\,
  \underbrace{\Big(\sum_{g=1}^{\Nf}(\rho_e^g)^p\Big)^{\frac{1}{p}-1}(\rho_e^f)^{p-1}}_{\text{aggregation}}\\
  &\cdot\,
  \underbrace{\frac{1}{N_\mathrm{ip}}\sum_{i=1}^{N_\mathrm{ip}}}_{\text{integration}}
  \Bigl(
  \underbrace{\sens{H}{d}}_{\text{transition}}\,
  \underbrace{\sens{d(\bx_i;\bolds^f)}{s_j^f}}_{\text{distance}}
  \Bigr),
\end{split}
  \label{eqn:gradient_chain}
\end{equation}
\noindent for $1 \leq j \leq \Nsf, 1 \leq f \leq \Nf$, with the factors given in \eqnref{eqn:obj_grad}, \eqnref{eqn:H_grad} and \eqnref{eqn:sens_d_p_x}, the latter depending on which case yields the shortest distance.

With a state-based objective $J(\bu(\brho(\bolds)))$, the function factor $\sens{J}{\widetilde{\rho}_e}$ is a global field influenced by support, load and potentially also far-away features. Even though the transition factor $\sens{H}{d}$ vanishes outside the transition zone, such that $\sens{J(\bu(\brho))}{\widetilde{\rho}_e}$ is collected only within the zone, the global information remains. With our stateless $\sens{J(\brho)}{\widetilde{\rho}_e}$ for tracking and reward, there are no global effects and we therefore use the asymmetric transition zone.

\section{Second-order formulation}
\label{sec:second_order}
The clearly most dominant optimizer in topology optimization is the Method of Moving Asymptotes (MMA)~\cite{svanberg1987method}, and also the most common approach to feature mapping is to use MMA as the optimizer. So it is clear that MMA can effectively solve feature-mapping problems. Nevertheless, we take here the chance to discuss and apply second-order optimization for feature mapping, as our specific objective functions allow for cheap-to-evaluate exact Hessian information.
\subsection{Background on the Hessian}
In standard density-based topology optimization, e.g., the common compliance minimization problem, we have a huge number of design variables which control a finite-element-based state problem. It is practically prohibitive to compute or even store the dense full Hessian matrix. This can be shown with the SIMP compliance gradient, already augmented by the residual
\begin{align*}
  \total{J_\mathrm{compl}}{\rho_e} = -\bu^\top \sens{\bK}{\rho_e} \bu + \blmbd^\top \left(\bK \, \bu - \boldf \right)
\end{align*}
and, following the standard adjoint-approach, the second derivative is given as
\begin{align}
 \frac{\mathrm{d}^2 J_\mathrm{compl}}{\mathrm{d}\rho_e\,\mathrm{d}\rho_{e'}} & = \blmbd_e^\top \sens{\bK}{\rho_{e'}} \bu - \bu^\top \scndcross{\bK}{\rho_e}{\rho_{e'}} \bu
  \label{eqn:compliance_simp_hessian}\\
  \bK\,\blmbd_e & = 2 \sens{\bK}{\rho_e}\bu. \nonumber
\end{align}
So for every element $e$ we need a specific right-hand side and get the corresponding adjoint solution $\blmbd_e$. See \cite{RojasLabanda2016} and \cite{FuKennedy2023} for a detailed discussion of the Hessian in the context of SIMP and approaches to actually use second-order information for specific formulations.

%The approximations MMA constructs for each function are separable, their analytical Hessians are diagonal only and cheap to compute and store. 
The optimizer IPOPT~\cite{Waechter2006}, far less established in the structural optimization community, internally solves Newton-type problems and usually does a limited-memory approximation (L-BFGS) of the Hessian. As such, SIMP problems can be solved reasonably using IPOPT. For appropriate problems, e.g., with a small number of design variables, IPOPT also takes exact full Hessian matrices. Because of this option we mainly use IPOPT in this manuscript to present a fair comparison between the first- and second-order formulations.

The Hessian $\bH$ collects the entries $H_{ij}^{fg}$ \eqnref{eqn:hessian}, which differentiate any gradient of variable $s_i^f$ in feature $f$ with respect to $s_j^g$, a variable of feature $g$.
\begin{equation}
\begin{gathered}
  H_{ij}^{fg} = \frac{\partial^2 J(\brho(\bolds))}{\partial s_i^f \, \partial s_j^g}
  = \sens{}{s_j^g}\,\sens{J(\brho(\bolds))}{s_i^f}, \\
  1 \le i \le \Nsf, \; 1 \le j \le \Nsg, \; 1 \le f,g \le \Nf
\end{gathered}
  \label{eqn:hessian}
\end{equation}
We can characterize $H_{ij}^{fg}$ as follows. $H_{ii}^{ff}$ is the \textit{curvature}; it expresses the variable-specific scaling -- important as the variables $q_x, r, \ldots$ have widely varying impact. %MMA does this within few iterations by adjusting the asymptotes. 
Then we have \textit{within-feature coupling} blocks $H_{ij}^{ff}$ for $i\neq j$, which support moving and rotating a feature, with its variables acting as a team. Finally we have \textit{cross-feature coupling} blocks $H_{ij}^{fg}$ for $f\neq g$, which support coordinating the movement of different features. Under the tracking objective \eqnref{eqn:tracking_obj}, two features competing for the same target region repel each other already through the overshoot of the $p$-norm aggregation; the cross-feature coupling adds the information that both features move at the same time. We have non-zero cross-feature coupling blocks only when the features overlap the same density elements $\rho_e$.   

\subsection{Sensitivity analysis}
\label{sec:second-order_analysis}

$\bH$ is symmetric and in our examples with at most 10 features and without additional fading variables $50 \times 50$ in size. The derivation of the Hessian is straightforward, purely technical but tedious to do manually and difficult to write and to read in compact form. For the derivative of the distance function \eqnref{eqn:dist} we have to consider the three cases for the closest part of the capsule. Users might use automatic differentiation or code generation. We restrict the presentation here to sample cases to improve readability and refer to \cite{Jung2026arXiv} for the full derivation and the implementation within openCFS (see code availability section).

The Hessian assembles from the same chain as the gradient, each term
carrying the second derivative of exactly one layer:
\begin{equation}
\begin{split}
  & \frac{\partial^2 J}{\partial s_i^f\,\partial s_j^g}
  = \sum_{e=1}^{\Ne}\Bigg[
    \underbrace{\scnd{J}{\widetilde{\rho}_e}\;
      \sens{\widetilde{\rho}_e}{\rho_e^f}\sens{\rho_e^f}{s_i^f}\;
      \sens{\widetilde{\rho}_e}{\rho_e^g}\sens{\rho_e^g}{s_j^g}}_{\text{function}} \\
  & + \underbrace{\sens{J}{\widetilde{\rho}_e}\,
      \frac{\partial^2\widetilde{\rho}_e}{\partial\rho_e^f\partial\rho_e^g}\,
      \sens{\rho_e^g}{s_j^g}\,
      \sens{\rho_e^f}{s_i^f}}_{\text{aggregation}}
    + \underbrace{\delta_{fg}\,\sens{J}{\widetilde{\rho}_e}\,
      \sens{\widetilde{\rho}_e}{\rho_e^f}\,
      \frac{\partial^2\rho_e^f}{\partial s_i^f\partial s_j^f}}_{\text{feature}}
  \Bigg],
\end{split}
  \label{eqn:hessian_terms}
\end{equation}
for all $1 \le i \le \Nsf, \; 1 \le j \le \Nsg, \;\; 1 \le f,g \le \Nf$.
The feature term acts within a single feature only; therefore $\delta_{fg} = 1$ for $f = g$ and 0 otherwise. With ${\partial^2 J_\text{reward}}/{\partial \widetilde{\rho}_e^2} = 0$, the function curvature vanishes for the reward problem, while ${\partial^2 J_\text{track}}/{\partial \widetilde{\rho}_e^2} = 2$. The aggregation term has cross-feature coupling, but only when features overlap the same element $e$ as
\begin{equation}
  \frac{\partial^2\widetilde{\rho}_e}{\partial\rho_e^f\partial\rho_e^g}
  = (p-1)\left[\delta_{fg}\,\frac{(\rho_e^f)^{p-2}}{\widetilde{\rho}_e^{\,p-1}}
    - \frac{(\rho_e^f\rho_e^g)^{p-1}}{\widetilde{\rho}_e^{\,2p-1}}\right].
  \label{eqn:agg_scnd_deriv}
\end{equation}
For the reward problem, where the function curvature vanishes, \eqnref{eqn:agg_scnd_deriv} is the only source of cross-feature coupling. It expresses the saturation of the aggregation: a feature gains little where another one already covers the element, and the Hessian tells the step how quickly this sets in.
We apply density integration \eqnref{eqn:rho_e_num_int} directly and, based on \eqnref{eqn:rho_sens}, we have
\begin{equation}
\begin{split}
  \frac{\partial^2\rho_e^f}{\partial s_i^f\partial s_j^f}
    = \frac{1}{N_\text{ip}}\sum_{k=1}^{N_\text{ip}} \Bigg[
     & \scnd{H(d(\bx_k; \bolds^f))}{d}\,\sens{d(\bx_k)}{s_i^f}\sens{d(\bx_k)}{s_j^f} \\
     & + \sens{H(d(\bx_k))}{d}\frac{\partial^2 d(\bx_k)}{\partial s_i^f\partial s_j^f}\Bigg].
\end{split}
 \label{eqn:rho_hess}
\end{equation}
The derivatives of the transition function, $\sens{H}{d}$ and $\scnd{H}{d}$, are given
by~\eqnref{eqn:H_grad} and~\eqnref{eqn:H_scnd}. The derivatives of the distance function follow trivially from the three-case signed distance, see~\eqnref{eqn:sens_d_p_x} for $\sens{d}{p_x}$. Each case has its own second derivatives. In the case of the shortest distance to the $P$-cap, we have
\begin{align*}
  \frac{\partial^2 d}{\partial p_x\,\partial p_x} = \dfrac{(y-p_y)^2}{\norm{\bx-P}^3}, \;
  \frac{\partial^2 d}{\partial p_x\,\partial p_y} = -\dfrac{(x-p_x)(y-p_y)}{\norm{\bx-P}^3}.
\end{align*}
\noindent See~\cite{Jung2026arXiv} for the full set of derivatives.
%The cap entries follow from the projector form $\frac{1}{R}\bigl(\bI-\hat{\bn}\hat{\bn}^{\T}\bigr)$, $\hat{\bn}=(\bx-P)/R$.

As neither objective is convex, the exact Hessian is in general indefinite. IPOPT handles this by its inertia correction, adding a multiple of the identity to the Hessian block of the KKT system whenever the factorization reports the wrong inertia; we never modify or convexify the exact Hessian ourselves. Second-order information changes the steps, not the landscape: the sensitivities still vanish outside the transition zone, and the local minima which follow from this remain.

\subsection{Generalization to state problems}
\label{sec:general_Hessian}
It is outside the scope of this work, but for completeness we also give here the extension for arbitrary state-based functions like the compliance. We then have to consider the general term 
\begin{equation*}
  \scndcross{J}{\rho_e}{\rho_{e'}},
\end{equation*}
which is generally a dense matrix and, in the compliance case, computed by the standard adjoint approach \eqnref{eqn:compliance_simp_hessian}, so that we need to solve for $\Ne$ right-hand sides.
The function block in \eqnref{eqn:hessian_terms} then becomes non-local: the diagonal
objective curvature $\scnd{J}{\widetilde{\rho}_e}$ is replaced by the dense mixed term
$\scndcross{J}{\widetilde{\rho}_e}{\widetilde{\rho}_{e'}}$ and a second element sum enters,
\begin{equation*}
   \frac{\partial^2 J}{\partial s_i^f\,\partial s_j^g} = \sum_{e=1}^{\Ne} \sum_{e'=1}^{\Ne}
    \scndcross{J}{\widetilde{\rho}_e}{\widetilde{\rho}_{e'}}\;\sens{\widetilde{\rho}_e}{s_i^f}\;
    \sens{\widetilde{\rho}_{e'}}{s_j^g} + \sum_{e=1}^{\Ne} \Bigg[ \ldots \Bigg],
  %\label{eqn:hessian_obj_general}
\end{equation*}
while the aggregation and feature terms of \eqnref{eqn:hessian_terms}, which only contract the first-order sensitivity $\sens{J}{\widetilde{\rho}_e}$, keep their single sum over $e$ and stay unchanged. For simplicity we assume a linear mapped density $\rho_e$ in the finite element solution, hence the second derivative of $\bK$ in \eqnref{eqn:compliance_simp_hessian} vanishes. With the sparse shape Jacobian (zero outside the transition zones)
\begin{equation*}
  D_{ei}^f = \sens{\widetilde{\rho}_e}{s_i^f}= \sens{\widetilde{\rho}_e}{\rho_e^f}\,\sens{\rho_e^f}{s_i^f} \in \mathbb{R}^{\Ne \times \Ns},
\end{equation*}
the adjoint solution $\blmbd_e = 2\,\bK^{-1}\sens{\bK}{\widetilde{\rho}_e}\,\bu$ of \eqnref{eqn:compliance_simp_hessian} and the symmetry of $\bK$, $\bK^{-1}$ and $\sens{\bK}{\widetilde{\rho}_e}$, the double sum becomes
\begin{equation}
\begin{split}
  & \sum_{e=1}^{\Ne}\sum_{e'=1}^{\Ne}
      \blmbd_e^\top \sens{\bK}{\widetilde{\rho}_{e'}}\,\bu \;\; D_{ei}^f\, D_{e'j}^g \\
  ={} & 2 \sum_{e=1}^{\Ne} D_{ei}^f\, \bu^\top \sens{\bK}{\widetilde{\rho}_e}\;
      \bK^{-1} \sum_{e'=1}^{\Ne} D_{e'j}^g\, \sens{\bK}{\widetilde{\rho}_{e'}}\,\bu \\
  ={} & 2\,(\ba_i^f)^\top \bg_j^g ,
\end{split}
  \label{eqn:hessian_obj_factorized}
\end{equation}
with the pseudo-load $\ba_i^f = \sum_{e=1}^{\Ne} D_{ei}^f\,\sens{\bK}{\widetilde{\rho}_e}\,\bu$ and the state derivative $\bK\,\bg_i^f = \ba_i^f$. Collected over all variables this is $2\,\bA^\top \bG \in \R^{\Ns \times \Ns}$, with $\bA = [\,\ba_i^f\,]$ and $\bG = [\,\bg_i^f\,] \in \R^{\Nu\times\Ns}$, one column per feature variable, and the global system matrix $\bK \in \R^{\Nu \times \Nu}$. This requires $\Ns$ solutions with the same system matrix $\bK$. With the common direct solvers for 2D, based on some form of triangular factorization, these are very cheap operations. The same factorization has been used for SIMP with a linear filter matrix in place of our shape Jacobian in \cite{RojasLabanda2016} and \cite{FuKennedy2023}, where the latter uses it to identify the positive semi-definite part of the Hessian rather than to evaluate it.

\section{Model simplification by feature consolidation}
\label{sec:geometry_variable}
In \cite{Norato2015}, the feature density before aggregation is modeled by a feature-specific \emph{size variable} $\alpha_f \in [0,1]$ as
\begin{equation*}
  \hat{\rho}_e^f = \alpha_f^q \, \rho_e^f(\bolds),
\end{equation*}
where $q > 1$ has a penalizing effect for compliance minimization, just as in SIMP optimization. This also implies that a volume constraint needs to be evaluated without penalization. This concurrently gives a result with a minimal number of features, as superfluous features are simply scaled out. %Without such a mechanism, unnecessary features usually shrink and/or overlap other features. However, we limit shrinking by a minimal feature profile with respect to the transition zone $r_\mathrm{min} > a$.

For the tracking formulation, we cannot use the same penalization approach. We also do not want to use a concave binarization function $\sum_f^{\Nf}  \alpha_f (1-\alpha_f)$. Instead, we propose consolidation as a post-tracking problem
\begin{equation}
\begin{aligned}
  \min_{\bz\,\in\, U_\mathrm{ad}} \quad
     & \sum_{f=1}^{\Nf} \alpha_f\\
  \text{s.t.}\quad
     \sum_{e=1}^\Ne \left( \rho_e^* - \widetilde{\rho}_e(\bolds, \bs{\alpha})\right)^2 & \leq (1 + \varepsilon) J_\mathrm{track}^* \\
      l^f(\bolds) &\ge \ell_\mathrm{min} \quad \forall\, 1\leq f \leq \Nf.
\end{aligned}
\label{eqn:consolidation_formulation}
\end{equation}
Following \cite{Norato2015}, we extend the design space by $\alpha_f$
\begin{equation*}
\begin{gathered}
\bz
=
\left(
(\bz^1)^\top,\dots,(\bz^{\Nf})^\top
\right)^\top \in \R^{\Ns + \Nf}, \\
\bz^f = (\bolds^f,\alpha_f)^\top
\quad \forall f=1,\dots,\Nf,
\end{gathered}
\end{equation*}
and set $U_\mathrm{ad}$ to the variable bounds of \eqnref{eqn:full_formulation} plus $0 \leq \alpha_f \leq 1$. The aggregated density \eqnref{eqn:pnorm} now becomes
\begin{equation}
  \widetilde{\rho}_e(\bz) = \left(\sum_{f=1}^{\Nf} \left(\alpha_f \,  \rho_e^f\right)^p\right)^{1/p}.
  \label{eqn:pnorm_alpha}
\end{equation}
Note that here and in the sum \eqnref{eqn:consolidation_formulation}, $\alpha_f$ enters linearly, in contrast to $\alpha_f^q$. $J_\mathrm{track}^*$ is the resulting tracking value of a previous optimization, the relaxation is problem dependent, e.g., 5 percent. Here we call $\alpha_f$ the \emph{fading variable} and, although every feature has its own variable, we do not call it a shape variable, as it has no impact on the distance. We do not binarize the fading variables. A feature may also settle at an intermediate value, as scaling it down can reduce the overshoot of the aggregation and thereby relax the tracking constraint, see \figref{fig:synth-consolidation}.

The sensitivity chains of \secref{sec:second-order_analysis} carry over to $\bz^f$ by a straightforward derivation. As $\hat{\rho}_e^f = \alpha_f\,\rho_e^f$, the shape gradients simply gain the factor $\alpha_f$, while for $\alpha_f$ the integration, transition and distance factors collapse to the feature density $\rho_e^f$ itself. The linear scaling keeps the second order cheap: the pure $\alpha_f$ curvature vanishes, the mixed one is first-order data, so no new second derivatives arise -- unlike with a penalized $\alpha_f^q$. The shape variable blocks are \eqnref{eqn:hessian_terms} with the function and aggregation terms scaled by $\alpha_f\,\alpha_g$ and the feature term by $\alpha_f$.

Note that a faded feature becomes degenerate: as $\alpha_f \to 0$, its contribution to $\widetilde{\rho}_e$ disappears and with it all derivatives of the tracking functional with respect to its shape variables, together with the associated Hessian blocks. Its shape variables are then undetermined at the solution, the Hessian is singular in these directions and the fading is practically irreversible, which is why we apply the consolidation only as a post-processing stage. 

\section{Staged optimization strategy}
\label{sec:stages}

\begin{figure}[tbp]
  \centering
  \includegraphics[width=0.9\linewidth]{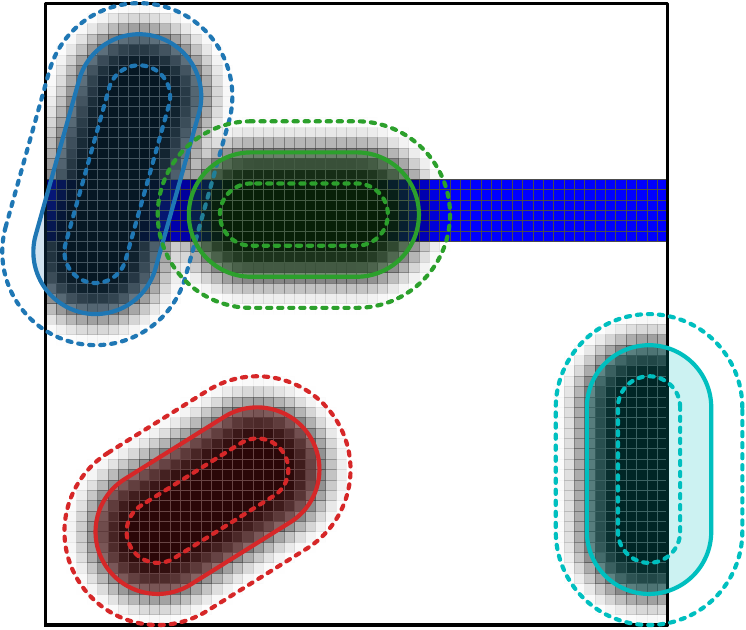}
  \caption{The blue horizontal bar is the target density $\brho^*$. In the initial configuration, the blue (b) and green (g) feature overlap with the target, the red (r) and cyan (c) features are far away, features (b) and (c) are only partially within the domain.}
  \label{fig:synth-setup}
\end{figure}

Our synthetic problem setup is depicted and explained in \figref{fig:synth-setup}. The domain is the unit square, discretized by 60 $\times$ 60 elements. The initial profile of the features is $r = 0.1$ and the transition zone is $a + a = 0.1$. The nominal width $2\,r$ is shown as full lines, the solid width as inner dashed lines, matching the width of the horizontal target bar. Numerical integration of $\rho_e$ is with $5 \times 5$ integration points. For all problems we require a minimal feature length of $0.01$ to prevent singularities in the derivatives.

As optimizers, we use for the first-order gradient case IPOPT \citep{Waechter2006} and GCMMA \citep{svanberg:2002:GCMMA}, the globally convergent variant of the standard Method of Moving Asymptotes, in the implementation of J\'er\'emie Dumas (\url{github.com/jdumas/mma}). For the second-order case we use again IPOPT but provide the exact Hessians.
\subsection*{Initial tracking (a failing attempt)}
\begin{figure}[tbp]
  \centering
  \includegraphics[width=\linewidth]{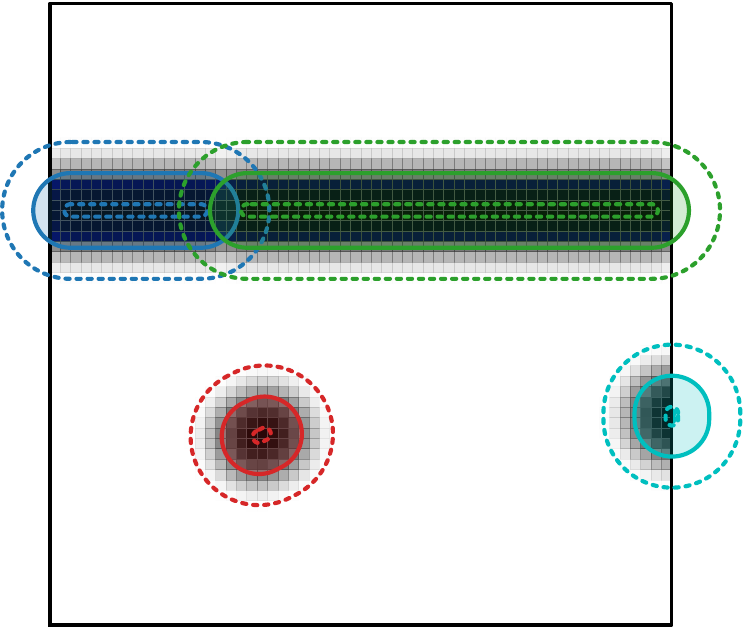}
  \caption{Tracking result from the initial configuration of \figref{fig:synth-setup}: Features (r) and (c) do not see the blue target and shrink / try to leave the domain.}
  \label{fig:synth-tracking}
\end{figure}

To motivate the usage of the reward function, we start with a study of the tracking objective \eqnref{eqn:tracking_obj} with a classical symmetric transition zone. We perform $p$-norm aggregation with a moderate $p=4$, hence overlapping is significantly larger than one and penalized by the tracking function. In \figref{fig:synth-tracking}, we see the expected result: the target is occupied by the features (b) and (g) in their initial order. They overlap only by their gray regions. The isolated features (r) and (c) do not see the target and more or less stay where they are and shrink to their minimal size to avoid the penalty of the tracking function. Feature (c) initially sees the domain boundary and also moves to the side. These results are more or less identical for the first-order case, the shown design is reached in fewer than 15 iterations. 

Surprisingly, IPOPT with the exact Hessian finds the target for all four features. This is likely due to a slight numerical drift of its internal barrier function for the constraints, and we consider this behavior to be incidental and therefore do not show it. A reference computation with the SciPy optimizer \texttt{trust-constr} also fails to bring features (r) and (c) to the target.

\begin{figure}[tbp]
  \centering
  \includegraphics[width=0.45\linewidth]{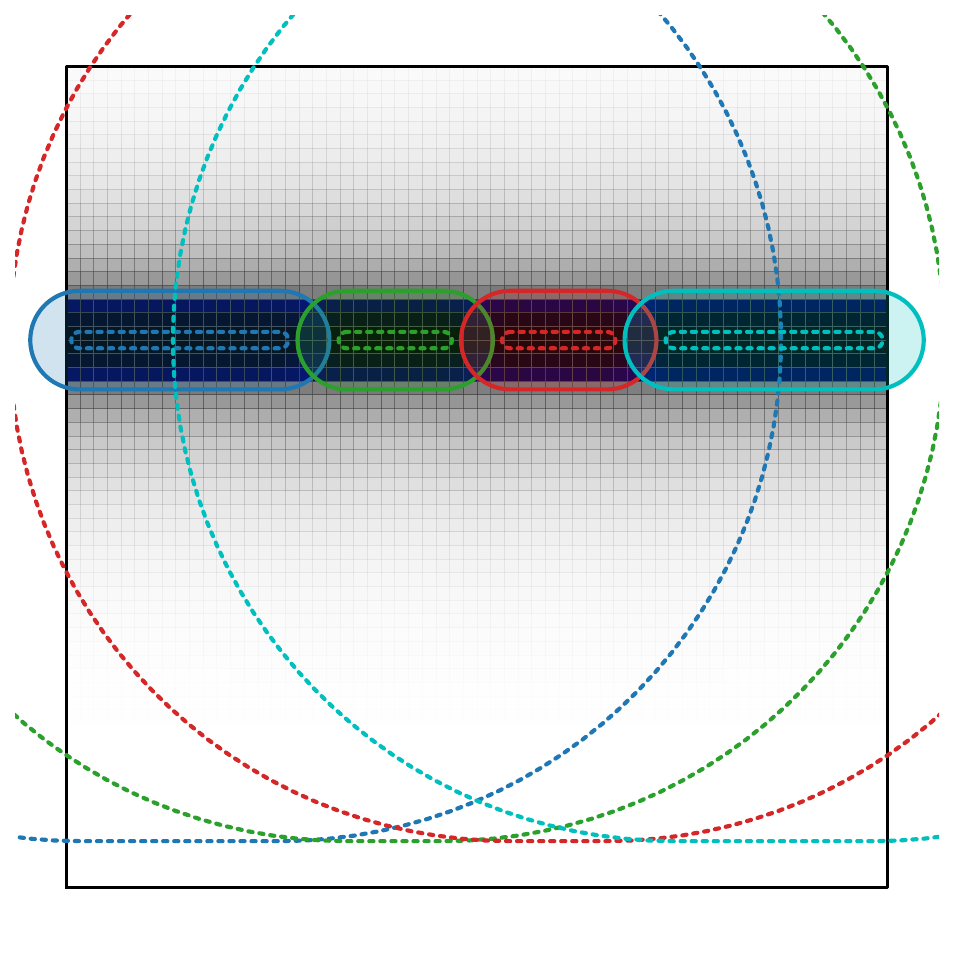} \,
  \includegraphics[width=0.45\linewidth]{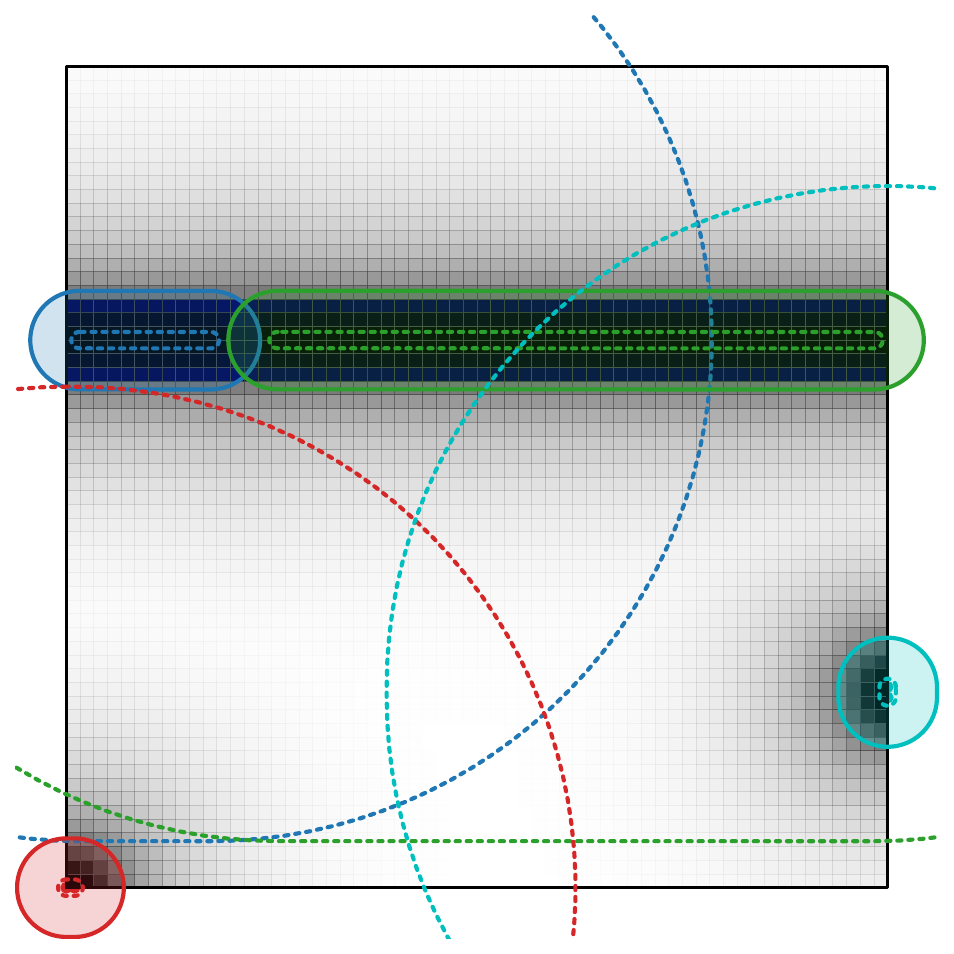} 
  \caption{Tracking with asymmetric transition zone, starting from the configuration in \figref{fig:synth-setup}. The left design is for the IPOPT Hessian case (d2). The IPOPT gradient (d1) result is similar. The right image shows the gradient case with GCMMA, which converges to a local minimum.}
  \label{fig:synth-asmytracking}
\end{figure}

Next, we run the same problems with an extension of $b=a+0.5$ and the quintic B\'ezier transition zone, see \figref{fig:synth-asmytracking}. Although all features see the target, GCMMA converges to a local minimum where the isolated features (r) and (c) are driven partially out of the domain. 

The value of the tracking function can be reduced by the intended overlapping with the target or by minimizing the area of the feature (shrinking). However, overlapping of features is also penalized \eqnref{eqn:p_norm_max}. This makes the tracking function, especially as initial stage, rather restrictive and sensitive to local minima.

\subsection*{Initial reward stage}
Our proposal to mitigate the tendency of the tracking function to get trapped in local optima is the reward function \eqnref{eqn:reward_obj}. The reward function gives a benefit when features overlap the target, but does not penalize non-overlapping. It also allows features to overlap each other: by the overshoot \eqnref{eqn:p_norm_max}, stacking features on the target even yields a slightly higher reward. We increase the parameter for the $p$-norm aggregation to $p=8$; however, the actual value is somewhat specific to the target geometry. Assume a more complex target than the trivial bar here, for instance one of the numerical examples in \secref{sec:results}. Let one part of the target with $n_A$ grid elements be already overlapped by $k$ features. Then for an additional feature, and ignoring gray elements, it is only beneficial to overlap a different blank target with $n_B$ elements if 
\begin{align}
  n_B > n_A \left( (k+1)^{1/p} - k^{1/p}\right).
  \label{eqn:n_B_ineq_reward}
\end{align} 

We furthermore fix the profile of the features to their initial value $r=0.1$. Enlarging to maximal features would be a trivial optimum.

\begin{figure}[tbp]
  \centering
  \begin{minipage}[c]{0.42\linewidth}
    \centering
    \includegraphics[width=\linewidth]{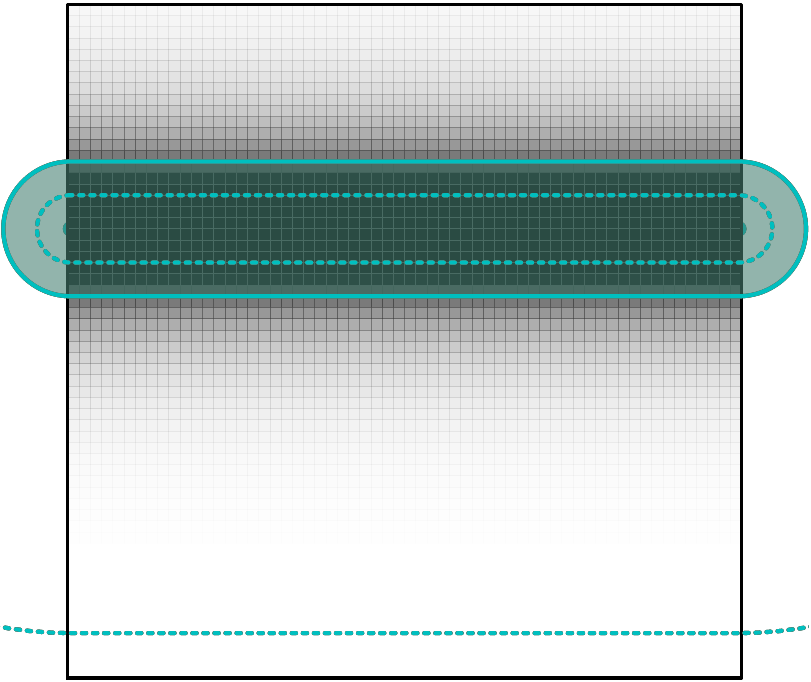}
  \end{minipage}\,%
  \begin{minipage}[c]{0.53\linewidth}
    \centering
    \includegraphics[width=\linewidth]{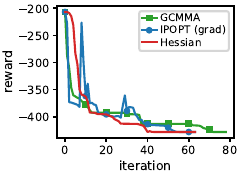}
  \end{minipage}
  \caption{Optimal solution for the reward function with asymmetric transition zone, starting with the design in \figref{fig:synth-setup}. All four features overlap with all three optimizers.}
  \label{fig:synth-reward}
\end{figure}

In \figref{fig:synth-reward}, we show the optimization result for the reward problem. All runs end in a similar design. For the Hessian case  $P = (0.004306, 0.666667)$ and for the GCMMA case $P = (0, 0.666667)$ for all features. For first-order IPOPT, we get $P_x \in [0.0191, 0.0231]$ and $P_y \in [0.6660, 0.6668]$ when converging with our very relaxed tolerance $0.01$. All reward values coincide to 6 digits, those of the Hessian and GCMMA case even to 9. The identical overlapping of all features is correct for the reward function, but it becomes an issue for GCMMA when overlapping is penalized by a subsequent tracking problem.

As can also be seen in the convergence plot, for the Hessian and GCMMA case, the features are roughly moved and oriented one after the other.

\subsection*{Bridging as second stage}
\begin{figure}[tbp]
  \centering
  \begin{subfigure}[b]{0.48\linewidth}
    \centering
    \includegraphics[width=\linewidth]{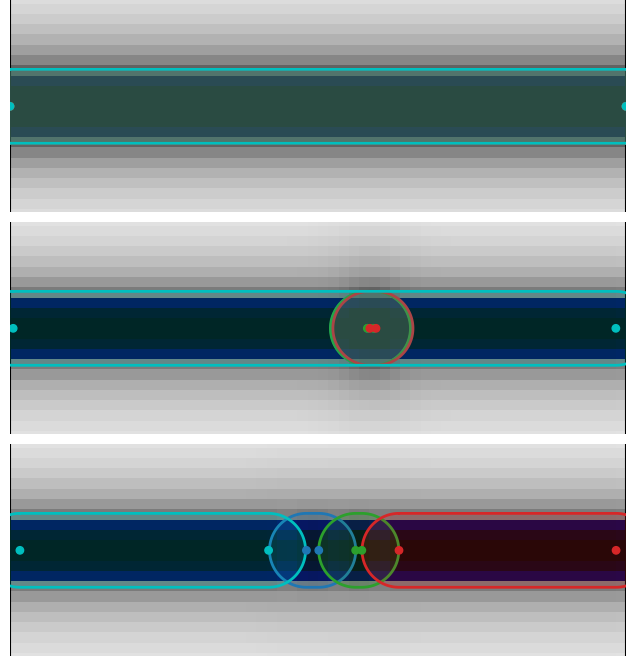}
    \caption{bridging (2nd stage)}
  \end{subfigure} \,
  \begin{subfigure}[b]{0.48\linewidth}
    \centering
    \includegraphics[width=\linewidth]{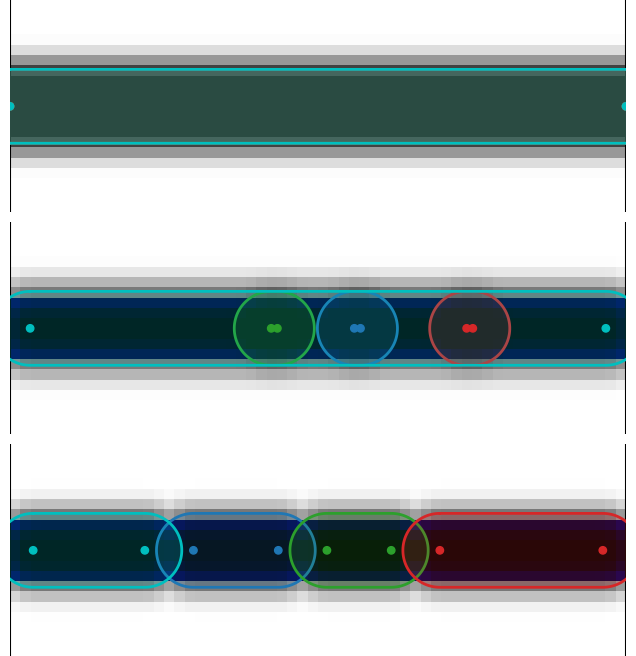}
    \caption{tracking (3rd stage)}
  \end{subfigure}
  \caption{From top down: GCMMA, IPOPT gradient, IPOPT Hessian. The overlapping of the non-Hessian solutions is not optimal with respect to the tracking objective.}
  \label{fig:synth-bridging}
\end{figure}

The reward stage has brought the features onto the target, but with fixed profiles and unpenalized overlap. The bridging stage therefore starts from the reward design and performs tracking, while we keep the asymmetric transition zone and the aggregation parameter $p=8$: removing the extension together with the switch to tracking would leave features which still have to travel to a non-overlapping area without a gradient. The results are shown in \figref{fig:synth-bridging}a. 

The benefit of the intermediate stage was observed in a random-start study in \citep{Jung2026arXiv}, which we do not repeat here.

%It depends on the intention how the bridging results are judged. We assume the tracking problem to prepare a design for a subsequent usage, e.g., compliance minimization. The tracking optimization is just an intermediate step here, and the function value likely has no real significance for the actual downstream application. All features are aligned with the target density, the profiles are accordingly adjusted.

The GCMMA start design was perfectly overlapped and the optimizer could not separate the features. What was beneficial with the reward function is here a clearly poor local saddle point with a tracking value of $131.51$. 

The IPOPT gradient problem also ends in a local optimum. One feature remains large to overlap the full target, the other three features shrink to their minimal size and stack. The tracking value is $100.97$. GCMMA finds a similar design when starting with the non-exact reward design. 

With the Hessian formulation IPOPT finds the best solution in the sense of the tracking objective ($99.26$) with features only overlapping by their transition zones. Due to the large transition zones, a smaller value cannot be achieved. A similar design is also found when starting with the other reward solutions.

\subsection*{Tracking as third stage}
Bridging already used the tracking function. We now remove the extension for the third and, in principle, final stage and set $p=4$ for the $p$-norm aggregation. The relaxed $p$-norm removes a little of the nonlinearity and increases the repelling effect for the features, but overlapping is still not aggressively prevented. As expected, the change for this trivial example is minimal, see \figref{fig:synth-bridging}b.

Removing the extension reduces the obtained objective value, while reducing $p$ from 8 to 4 penalizes remaining overlapping more. As GCMMA cannot remove the exact overlap, its obtained $J_\mathrm{track}^*$ remains large at 120.7. Gradient IPOPT and Hessian IPOPT obtain 54.1 and 44.3, respectively. Because of our rather large symmetric transition zone, we cannot expect zero tracking.
%and the huge GCMMA value can be explained by the overlapping of four features $4^{1/4} \approx 1.41$.

\subsection*{Consolidation as optional fourth stage}
We adopt Norato's idea of scaling/fading variables $\alpha_f$, see \secref{sec:geometry_variable}; however, this is an extension of standard feature-mapping formulations, and we apply our consolidating feature minimization formulation \eqnref{eqn:consolidation_formulation} only as an option. 

\begin{figure}[tbp]
  \centering
  % the design is raised such that its frame is centered on the axes box of the
  % plot (the plot image carries the x-label below its axes box)
  \raisebox{13.2pt}{\includegraphics[width=0.4\linewidth]{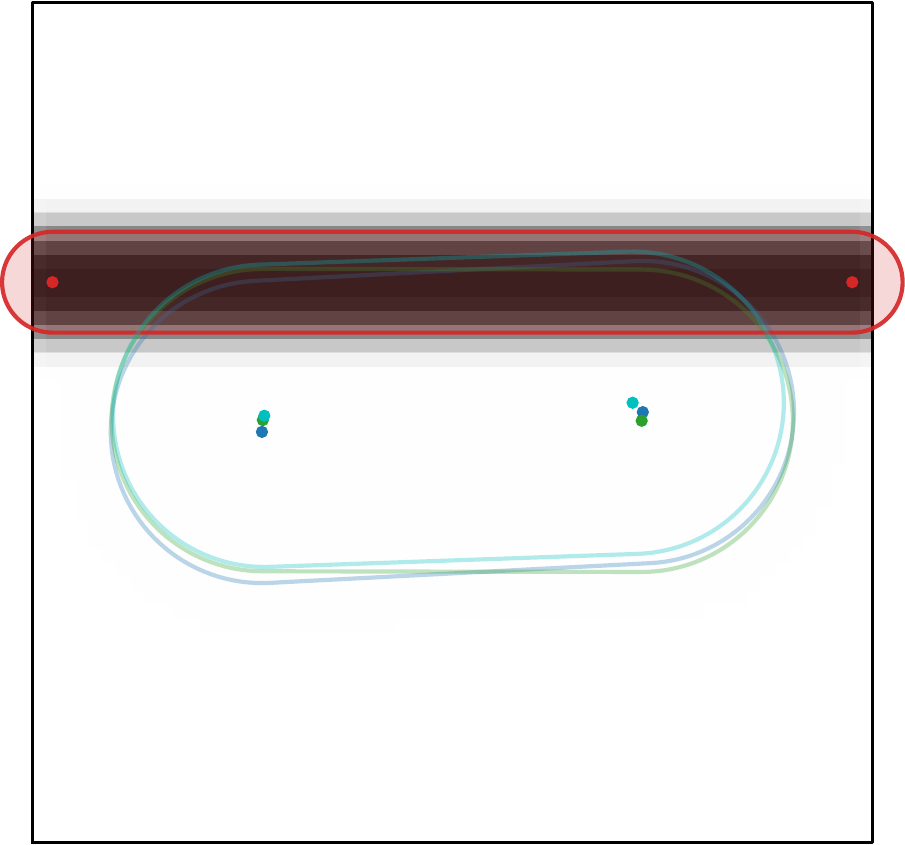}}\,%
  \includegraphics[width=0.55\linewidth]{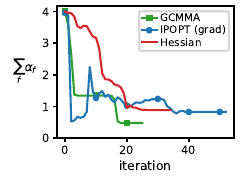}
  \caption{Optional consolidation stage: A typical solution of minimizing features via the fading variables $\alpha_f$ \eqnref{eqn:consolidation_formulation} (here the Hessian case). Three $\alpha_f$ become zero and one scales according to the individual tracking bound $J_\mathrm{track}^*$; its large value for GCMMA explains the low sum.}
  \label{fig:synth-consolidation}
\end{figure}

The $\varepsilon$ of the $(1 + \varepsilon)J_\mathrm{track}^*$ tracking bound with the $\sum_f \alpha_f$ objective is very much problem dependent and reflects how much the features still have to be rearranged. A larger relaxation helps. When the final fading variables $\alpha_f$ were not already close to binary 0 or 1, even a trivial threshold of 0.5 was sufficient in our numerical experiments.
A further option, outside the scope of this work, is to solve the actual optimization problem, e.g., compliance, in the sense of geometry projection with an $\alpha$ formulation, starting from the tracking or bridging design.
For all optimizers, the consolidation stage finds a single bar as a solution, even for the GCMMA case, see \figref{fig:synth-consolidation}. As the tracking bound for GCMMA is rather large due to the poor local optimum, three features become zero and the remaining feature can be intermediate. In that case the bound could even be tightened below $J_\mathrm{track}^*$, but a robust threshold gives the same optimal design, a single bar.

\section{Numerical results}
\label{sec:results}

We now apply the staged strategy of \secref{sec:stages} to the SIMP solutions of two linear elastic compliance benchmark problems: a single-load cantilever on a $60 \times 60$ mesh of the unit square, and a two-load five-bar mechanism on a $120 \times 60$ mesh of a 2 $\times$ 1 domain; the edge size is $h=0.01\overline{6}$. The SIMP problems are standard, with penalization exponent 3, volume fraction 0.5 and a density filter of radius $1.8\,h$ (cantilever) and $1.4\,h$ (five-bar). 

We use a simple initial feature design with five and ten features, respectively, whose asymmetric transition zones overlap (see \figref{fig:cant-pipeline}a and \figref{fig:fb-pipeline}a). %Our configuration is also known as quincunx
Other common initial configurations in the literature are cross-braced patterns or a regular lattice, in all cases with some overlap between neighboring features. Here we have no support or load which we want to cover initially by features.

We perform each pipeline with our exact Hessian implementation via IPOPT \citep{Waechter2006} and also first-order with the default L-BFGS mode of IPOPT and a globally convergent MMA \citep{svanberg:2002:GCMMA}. The shown designs are from the Hessian case; we only comment on differences from the other two optimizers. 

We generally use numerical integration with 5 $\times$ 5 integration points. The symmetric transition zone is 0.1, hence $a=b=0.05$. The standard asymmetric extension for the reward and bridging stage is 0.2, hence $b=0.25$. See the exception for the five-bar bridging stage below. With respect to the $p$-norm aggregation, we use $p=8$ for the reward stage to not reward overlapping too much. For the bridging stage (tracking function with asymmetric transition zone), we use $p=8$ or $p=4$, as the penalization of overlapping (drive for repelling) needs to be adjusted for the specific problem. For the symmetric tracking stage, $p=4$ is used. All settings apply equally to the three optimizers.

In all cases we have a minimal feature length of 0.05. When the profile is free, its bounds are from 0.06 to 0.5 and larger than the inner transition zone. 
\subsection{Single-load cantilever}
\label{sec:results_cantilever}
\begin{figure}[tbp]
  \centering
  %% the sub-captions may slightly overhang their box, so that
  %% "bridging/tracking" stays on a single line
  \captionsetup[sub]{margin=-4pt}
  \begin{subfigure}[b]{0.32\linewidth}
    \centering\includegraphics[width=\linewidth]{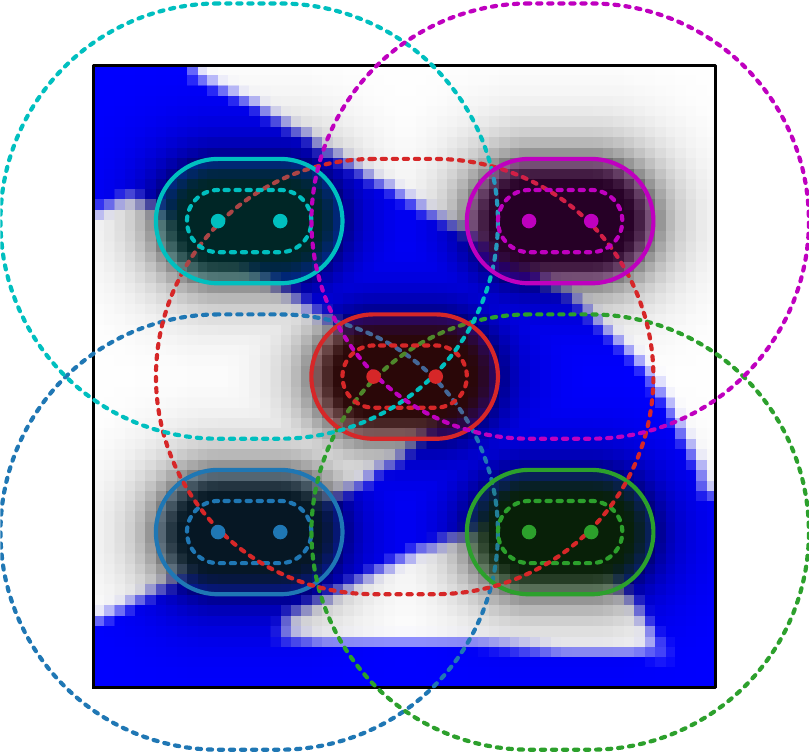}
    \caption{initial}
  \end{subfigure}\hfill
  \begin{subfigure}[b]{0.32\linewidth}
    \centering\includegraphics[width=\linewidth]{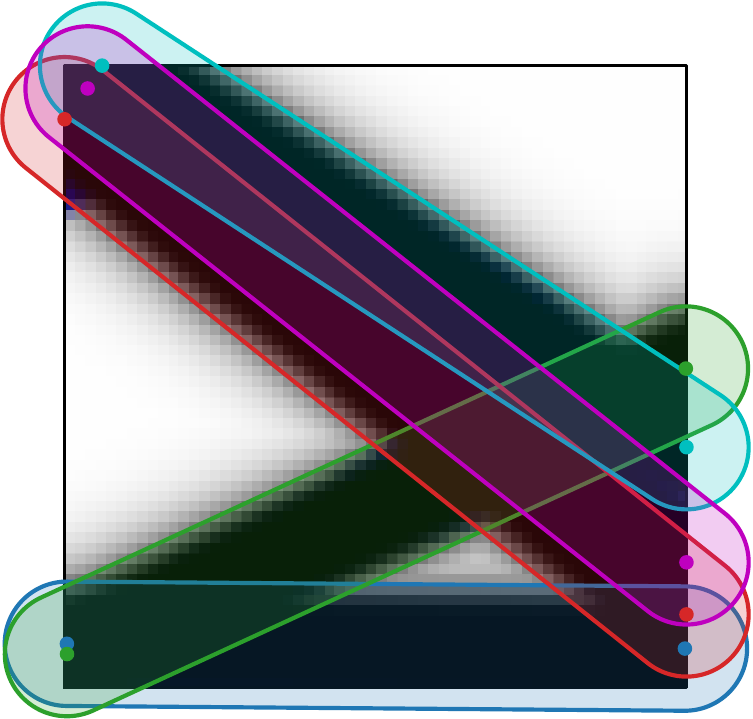}
    \caption{reward}
  \end{subfigure}\hfill
  \begin{subfigure}[b]{0.32\linewidth}
    \centering\includegraphics[width=\linewidth]{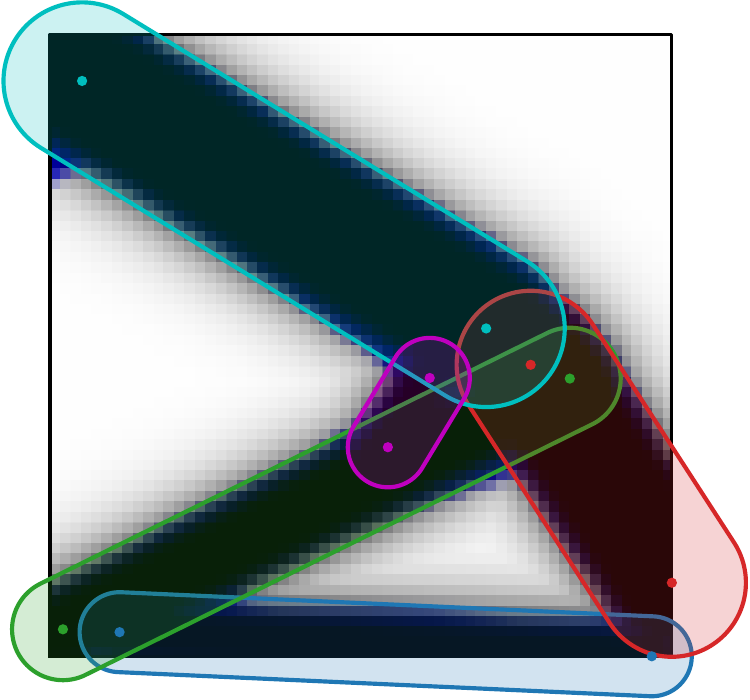}
    \caption{bridging/tracking}
  \end{subfigure}
  \caption{Cantilever pipeline using Hessian IPOPT. The bridging (shown) and tracking designs are basically identical.}
  \label{fig:cant-pipeline}
\end{figure}

In \figref{fig:cant-pipeline}, we show the staged pipeline for a standard single-load cantilever SIMP solution as target density. The design can obviously be represented by three, or better, four features, hence one of our five features should be superfluous.

The reward stage converges for all three optimizers in roughly 35 iterations to a design similar to \figref{fig:cant-pipeline}b. As overlapping is rewarded and covering void is not penalized, overlapping the slightly kinked diagonal with three features (with fixed profile) is the optimal solution.

For the bridging stage, the three optimizers show different behavior. The Hessian case in \figref{fig:cant-pipeline}c shows the best design; it converged in fewer than 50 iterations. GCMMA takes a little longer and cannot separate the single nearly superfluous feature but has two long features overlapping the diagonal. Note that we have no upper length bound for the features. IPOPT with first-order gradients tends to take large steps within the first iterations, hence the staging pipeline starting from the previous designs loses some of its motivation. With the given parameters, the lower horizontal bar is not covered, with some other $p$-aggregations it is, but this is not predictable. 

The tracking stage just switches to a symmetric transition zone and increases the repelling effect by lowering $p$ to 4. As expected, the designs show almost no visible change and convergence is therefore very fast within about 20 iterations (with the exception of L-BFGS IPOPT, which does not converge).

\subsection{Five-bar two-load}
\label{sec:results_fivebar}
\begin{figure*}[tbp]
  \centering
  \begin{subfigure}[b]{0.48\linewidth}
    \centering\raisebox{1ex}{\includegraphics[width=\linewidth]{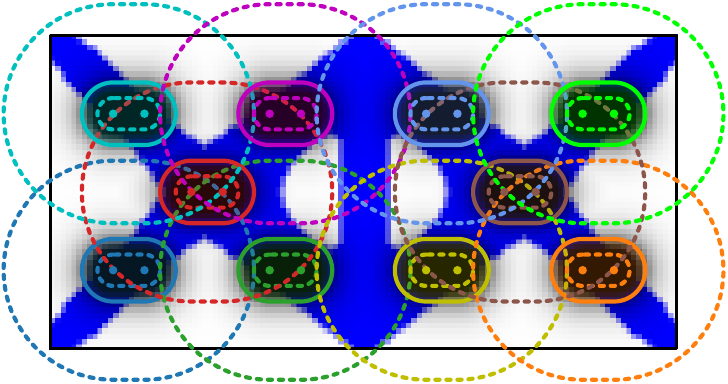}}
    \caption{initial}
  \end{subfigure} \,
  \begin{subfigure}[b]{0.48\linewidth}
    \centering\includegraphics[width=\linewidth]{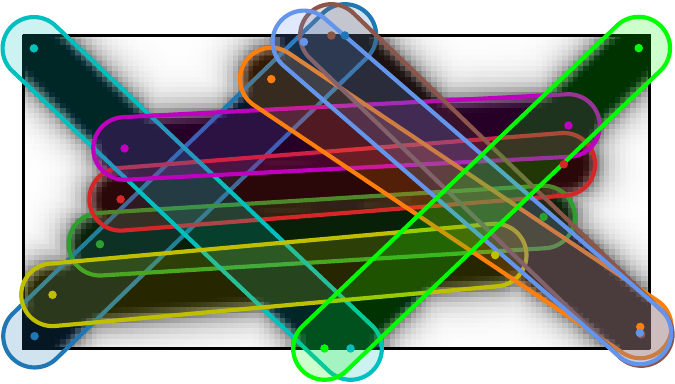}
    \caption{reward}
  \end{subfigure} \\
  \begin{subfigure}[b]{0.48\linewidth}
    \centering\includegraphics[width=\linewidth]{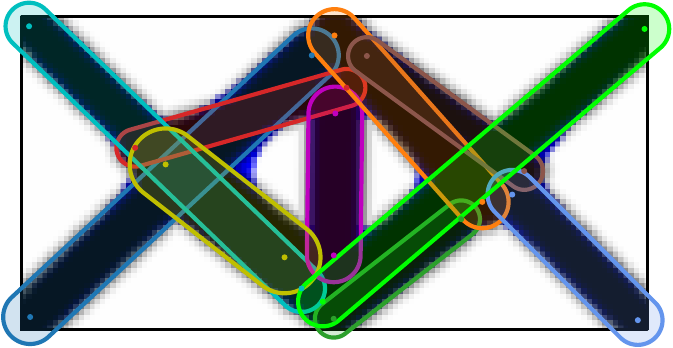}
    \caption{bridging/tracking}
  \end{subfigure} \,
  \begin{subfigure}[b]{0.48\linewidth}
    \centering\includegraphics[width=\linewidth]{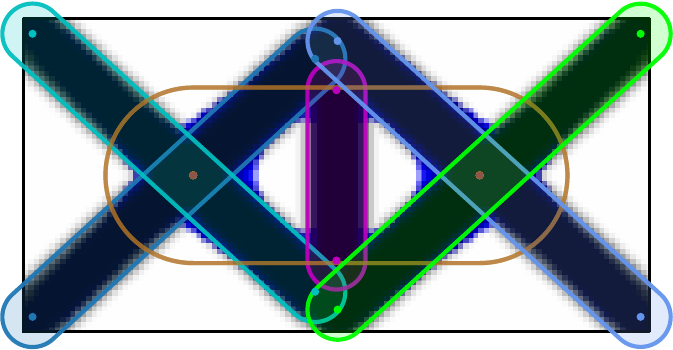}
    \caption{consolidation}
  \end{subfigure}
  \caption{Five-bar pipeline using Hessian IPOPT. The maximum length of the 10 features is limited to $\sqrt{2}$ to allow the diagonal bars but prevent covering the full domain. The designs in the bridging and tracking (shown) stages are basically identical. The optional consolidation formulation \eqnref{eqn:consolidation_formulation} reduces to five bars.}
  \label{fig:fb-pipeline}
\end{figure*}

The five-bar benchmark in \figref{fig:fb-pipeline} is a common two-load problem with full support at the corners and separate vertical loads at the centers of the long edges. We initialize with 10 features, so there are many potentially superfluous features. The problem is surprisingly difficult, as an optimal reward solution is far away from an optimal tracking solution.

For the reward problem we limit the length of the features to $\approx \sqrt{2}$, thus preventing the full domain from being covered by 10 horizontal features while still allowing for the diagonal bars of the target. The fixed profile roughly matches the bar width of the target. \figref{fig:fb-pipeline}b shows that the two crosses are nicely matched by four features and most of the remaining features are horizontal to fully capture the target while maximizing the overlap. Due to the maximal length restriction and symmetry of the design, the shown optimal design for the Hessian formulation is obviously not unique. Several features cover the same part of the target and can be shifted along it without changing the reward, so the objective is flat in these directions and the exact Hessian is singular there. This is a plausible reason why the second-order run takes the most iterations to settle, an effect which does not appear for the cantilever. All three optimizers basically find a similar design. L-BFGS IPOPT converges after 40 iterations, GCMMA after 80 and the Hessian formulation only after 120. 

For the bridging stage it is difficult to track the vertical center bar as it is covered by several vertical features at the end of the reward stage. The Hessian formulation succeeds, but with a slight increase of the extension to 0.25, hence $b=0.3$, all three optimizers find a solution basically similar to \figref{fig:fb-pipeline}c. L-BFGS IPOPT takes about 60 iterations to find the design but does not yet reach its convergence criteria. The Hessian formulation converges within 80 iterations and GCMMA finds the final design only after 150 iterations. 

Again, for the tracking stage, no significant visual design change can be observed and all optimizers converge within 25 iterations.

% ==========================================================================
\section{Conclusion}\label{sec:conclusion}
% ==========================================================================
The presented approach poses the search for a feature-mapping design aligned to a given density field as an optimization problem. The purpose of this task might be to use a SIMP design to generate an initial configuration for a successive feature-mapping optimization or to generally interpret given density fields by a given set of feature geometries. The proposed tracking function is a straightforward approach, which is, however, highly prone to local minima and therefore itself depends on a good initial feature configuration. Our remedy is an initial reward step, which promotes aligning features to the target density but without the penalization of void covering and feature stacking. Altogether, we propose three stages, starting with the asymmetric reward, followed by asymmetric tracking as a bridging stage and finally symmetric tracking, optionally extended by a consolidation stage based on the feature fading variable of the geometry projection method. Whether all three stages are really required depends very much on the use case. In the work underlying this paper \citep{Jung2026arXiv}, large problem sets with randomly generated target designs were also studied -- there the three stages proved to be a robust approach.

Tracking and reward are stateless, which reduces the sensitivity information strictly to the intermediate mapped pseudo density within the transition zone. We therefore propose an asymmetric transition zone based on B\'ezier curves together with an algorithm to automatically determine the B\'ezier parameters from the size of the transition zone. We provide a simple web tool to experiment with the B\'ezier curves at \url{https://am-ko.mi.uni-erlangen.de/bezier.html}.

The statelessness which led to the asymmetric transition zone pays off in the second order: the exact Hessians of the tracking, reward and consolidation functionals require no additional system solution and are cheap to evaluate. The Hessian for state-based feature mapping requires additional solutions of the state problem for each feature variable, which is usually very cheap in conjunction with direct solvers. For completeness, we give the formulation for the compliance case.

We compute the numerical examples with IPOPT and GCMMA in the first-order formulation and again IPOPT for the second-order formulation. In our examples, the latter is usually remarkably more robust with respect to parametrization; however, the gain is not dramatic. A precise study and comparison of the first- and second-order formulations is beyond the scope of this work. Such a comparison depends very much on the selection of optimizers (including their handling of indefinite Hessians) and on their stopping criteria. It should also include state-based problems like compliance, force inverters and possibly harmonic formulations.

\backmatter

\bmhead{Acknowledgements}
The initial part of this work goes back to the project FIONA (LuFo IV-1, FKZ: 20W1913F) funded by the German Federal Ministry for Economic Affairs and Climate Action (BMWK). Within this project, Jannis Greifenstein implemented and successfully tested the tracking problem.

\section*{Statements and Declarations}
\textbf{Funding } German Federal Ministry for Economic Affairs and
Climate Action (BMWK) under the project FIONA (LuFo IV-1, FKZ: 20W1913F)

\bigskip
\noindent \textbf{Conflict of interest } There are no conflicts of interest to declare that are relevant to the content of this article.

\bigskip
\noindent \textbf{Author contributions } The original tracking formulation was proposed by Michael Stingl. Patrick Jung did an extensive study based on his own implementation in his Master's thesis \textit{Optimizing Initial Feature-Mapping Variables from Given Designs via Tracking} \citep{Jung2026arXiv}, supervised by Fabian Wein and Michael Stingl. The thesis contains the three stages, the asymmetric transition zone and the Hessian formulation. Fabian Wein is the main author of the present work and contributed the B\'ezier-based asymmetric transition zones and the present Hessian formulations. Arash Moradian provided most of the numerical examples. All authors contributed to the editing and proofreading of the manuscript.

\section*{Replication of results}
Exact equations are included in the paper, making it straightforward for readers to replicate results. The methods are all integrated into the open source software openCFS \citep{mayrhofer2026opencfs}, see \url{https://opencfs.org} for links to binaries and
the source code. Representative test cases are given at \url{https://gitlab.com/openCFS/cfs/-/tree/master/Testsuite/TESTSUIT/Optimization/FeatureMapping}.
An interactive tool to explore the B\'ezier-based transition function is available at \url{https://am-ko.mi.uni-erlangen.de/bezier.html}.
The datasets generated and analyzed during the current study are available from the corresponding author on request.

\bibliography{tracking}

\end{document}